\def\Gn{\,\rule[-1mm]{2pt}{5mm}\,}

\def\BE{\begin{equation}}
\def\EE#1{\label{#1}\end{equation}}
\def\C{{\mathbb C}^3}
\def\R{{\mathbb R}^3}
\def\T{{\mathbb T}^3}
\def\e{{\rm e}}
\def\i{{\rm i}}
\def\d{{\rm d}}
\def\D{\mathfrak{D}}
\def\wb{\widetilde{\bf b}}
\def\wv{\widetilde{\bf v}}
\def\hB{\widehat{\bf B}}
\def\hV{\widehat{\bf V}}
\def\O{\mathscr{O}}
\def\P{\mathscr{P}_{\bf n}}
\def\S{\mathfrak{S}}
\def\ba{\begin{align*}}
\def\be{\begin{align}}
\def\ks{\bm\xi}
\def\ha{\widehat{\bf a}}
\def\hu{\widehat{\bf u}}
\def\rf#1{(\ref{#1})}
\def\se#1{\begin{subequations}\label{#1}
\renewcommand{\theequation}{\theparentequation.\arabic{equation}}}
\documentclass[a4paper,12pt]{Melsarticle}
\usepackage{mathtools,amssymb,textcomp,wasysym,bm,mathrsfs,float}
\begin{document}
\oddsidemargin -1cm
\evensidemargin -1cm
\begin{center}
{\bf Space analyticity and bounds for derivatives of solutions
to the evolutionary equations of diffusive
magnetohydrodynamics}\footnote{%Zheligovsky V. Space analyticity and bounds for
%derivatives of solutions to the evolutionary equations of diffusive
%magnetohydrodynamics.
Mathematics 9, 1789 (2021). https://doi.org/10.3390/math9151789}

~

Vladislav~Zheligovsky

~

{\it Institute of Earthquake Prediction Theory and
Mathematical Geophysics, Russian Ac. Sci.,\\
84/32 Profsoyuznaya St, 117997 Moscow, Russian Federation}\\
\rule{\textwidth}{1pt}\vspace*{-1.25\baselineskip}\end{center}
\section*{Abstract}
In 1981, Foias, Guillop\'e and Temam proved a priori estimates for
arbitrary-order space derivatives of solutions to the Navier--Stokes equation.
Such bounds are instructive in the numerical investigation of intermittency often
observed in simulations, e.g., numerical study of vorticity moments
by Donzis et al. (2013) revealed depletion of nonlinearity that may be
responsible for smoothness of solutions to the Navier--Stokes equation.
We employ an original method to derive analogous estimates for space derivatives
of three-dimensional space-periodic weak solutions to the evolutionary
equations of diffusive magnetohydrodynamics. Construction relies
on space analyticity of the solutions at almost all times. An auxiliary
problem is introduced, and a Sobolev norm of its solutions bounds from below
the size in $\C$ of the region of space analyticity of the solutions
to the original problem. We recover the exponents obtained earlier
for the hydrodynamic problem. The same approach is also followed here to derive
and prove similar a priori bounds for arbitrary-order space derivatives
of the first-order time derivative of the weak MHD solutions.

\noindent\rule[.5\baselineskip]{\textwidth}{1pt}\vspace*{-\baselineskip}

\section{Introduction}\label{intr}

A standing problem of the analytical study of turbulence is to derive from
the basic equations of hydrodynamics, the~Euler and Navier--Stokes equations,
the empirical relations characterising this phenomenon. This requires
a profound understanding of the behaviour of small-scale structures in flows,
which is also necessary to achieve progress in pure mathematical problems
such as to identify the class of functions, in which existence and uniqueness
of solutions is guaranteed, or to answer the related question whether
singularities can develop at a finite time in the solutions.

A possible approach to addressing these problems consists
of obtaining information on norms of high-order derivatives of the solutions:
the higher the order, the more the respective norms are controlled by
the small-scale components of the solutions. The energy inequality
\BE{1\over2}\|{\bf V}\|_0^2+\nu\int_0^T\|{\bf V}\|_1^2\,\d t\le
\|{\bf V}^{\rm(init)}\|_0^2\EE{NSe}
for solenoidal solutions to the Navier--Stokes equation
\BE{\partial{\bf V}\over\partial t}=\nu\nabla^2{\bf V}
-({\bf V}\cdot\nabla){\bf V}-\nabla P\EE{NS}
bounds the Lebesgue space $L^2(\Omega)$ norms of an incompressible fluid flow
$\bf V$ and its
spatial gradient only, and not of the second derivatives describing the action
of diffusivity. This led J.~Leray \cite{Le} and E.~Hopf \cite{Ho} (see also
\cite{La,Te,RRS,Ro}) to formulate the concept of weak solutions to the Navier--Stokes
equation -- namely, vector fields satisfying integral relations that are
obtained by scalar multiplying \rf{NS} by a sufficiently smooth solenoidal
test function with a finite support, integrating over the fluid volume $\Omega$
and transferring differentiation from the unknown solution $\bf V$ to the test
function by integration by parts. If the resultant
integral identity holds for all such test functions and $\bf V$ is sufficiently
smooth, it is simple to show that it also solves \rf{NSe}; such solutions
are called strong. Later it was shown \cite{GL,So} that second-order spatial
derivatives and the time derivative of a three-dimensional weak solution
to \rf{NS} as well as the gradient of pressure, that are involved
in \rf{NS}, do exist and belong to the Lebesgue space
$L^{5/4}([0,T],L^{5/4}(\Omega))$; the proof relies on the observation that
for a vector field obeying the energy bound \rf{NSe} the nonlinear term in
\rf{NS} belongs to this space. Existence of weak solutions was demonstrated
in \cite{Le,Ho}; existence of strong three-dimensional solutions is an open
question. While for an incompressible fluid residing in a bounded
domain~$\Omega$ uniqueness was proven for three-dimensional flows satisfying
suitable boundary conditions and belonging to the Lebesgue space
$L^p([0,T],L^q(\Omega))$, for which the Ladyzhenskaya--Prodi--Serrin condition
$2/p+3/q\le1$ \cite{KL,Pr,Se,ESS} holds, the energy bound \rf{NSe} for weak
solutions implies only $2/p+3/q=3/2$.

Due to importance of these mathematical questions, numerous papers were devoted
to the investigation of smoothness and spatial analyticity of solutions
to the Navier--Stokes equations. In the seminal work \cite{FT}, C.~Foias and
R.~Temam examined Gevrey class regularity of space-periodic solutions and
proved that three-dimensional flows, which initially have spatial gradients
in $L^2(\T)$, instantaneously become space-analytic, and, for a finite time,
the size of the region of analyticity in $\C$ is proportional to time
(a similar derivation in \cite{DT} serves for estimating the minimum length
scales in the flow and Fourier spectrum decay in terms of the instantaneous
rate of the bulk energy dissipation; see also \cite{DG,FMRT}).
Space analyticity persists while the $L^2(\T)$ norm of $\nabla\bf V$ remains
finite (for weak solutions this can be guaranteed for finite times only).

In the celebrated paper \cite{FGT}, C.~Foias, C.~Guillop\'e and R.~Temam
established another regularising effect of the Navier--Stokes equation,
manifested by new a priori estimates: for initial conditions of a minimum
regularity, the weak solutions admit the bounds
\BE\int_0^T\!\|{\bf V}\|_m^{\alpha_m}\,\d t<\infty\qquad\mbox{for~~}
\alpha_m={2\over2m-1}.\EE{NSb}
(Here and in what follows, $\|\cdot\|_m$ denotes the norm in the Sobolev space
$H_m(\T)$; it is essentially equivalent to the sum of the $L_2(\T)$ norms
of all derivatives of order $m$.) This result was derived in \cite{G1,G2}
by a different method relying on the so-called ladder inequalities, employed
for estimating the ``natural'' length scale developing in a forced flow
\cite{BDGM,G1,G2,G3,G4,DG}. Recent developments in the study of analyticity of
solutions to the Navier--Stokes equations are described in \cite{BGK}.
An ordinary differential equation (ODE) is studied in \cite{BF}, that governs
the evolution of the size in $\C$ of the region
of analyticity of the solution and involves the Gevrey class norms; this
is reminiscent of the approach \cite{Gev} that we follow here.
A bound from below for the size of the region of analyticity that vanishes
on a measure zero time set was constructed in~\cite{Ch}.

An important problem is to characterise the singularities presumably
developing in solutions to the equations of hydrodynamics and
magnetohydrodynamics (MHD). Citing \cite{Te},
``It was Leray's conjecture on turbulence, which is not yet proved nor
disproved, that the solutions to the Navier--Stokes equations do develop
singularities ... It seems useful to study the properties of weak solutions
of Navier--Stokes equations with the hope of either proving that they are
regular, or studying the nature of their singularities if they are not. ...
Of course'' the results ``would lose all of their interest if the existence
of strong solutions were demonstrated.'' J.~Leray \cite{Le} showed
that for any weak solution of the force-free Navier--Stokes equation there
exists at most a countable set of disjoint open time intervals $J_q$ such that
$J_0$ is infinite, $J_q$ are finite for $q>0$,
$\sum_{\,q>0}\sqrt{{\rm length}(J_q)}<\infty$, the Lebesgue measure
of the complement $[0,\infty)\backslash\cup_qJ_q$ is zero and the solution is
smooth in all space-time regions $J_q\times\R$. If a body force acts
on the fluid, the singularity set has the same structure \cite{FGT} (except for
the inequality on the lengths of the time intervals of smoothness does not
necessarily remain valid). Investigation of the partial regularity
of solutions to the Navier--Stokes equations was continued by V.~Scheffer
\cite{Sch1,Sch2,Sch3,Sch4} and culminated in the work
by L.~Caffarelli, R.~Kohn and L.~Nirenberg \cite{CKN}, who proved that
for any suitable weak solution of the Navier--Stokes equation on an open set
in space-time, the singular set has a zero Hausdorff measure $\mathscr{H}^1$.

Proven bounds are instructive in numerical analysis of the nature
of intermittency observed in solutions to hydrodynamic or MHD equations.
For instance, the numerical study \cite{DV} of vorticity moments
of solutions to the Navier--Stokes equations revealed depletion of nonlinearity
that may be responsible for smoothness of the solutions under investigation.

Existence of weak solutions to equations of diffusive magnetohydrodynamics
was proven in~\cite{DL}. The large-time behavior of a solution
to the Navier--Stokes equation \cite{FTn} or an MHD solution \cite{ST} is
completely determined, if it is known in a sufficiently large, but finite set
of points in the fluid region. Since the nature of the quadratic nonlinearities
in the magnetic induction and
Navier--Stokes equations (in the MHD case, the latter involving
the Lorentz force acting on the electrically conducting fluid) is the same,
most results for the hydrodynamic
Navier--Stokes equation can be generalised, often straightforwardly,
to encompass the system of equations of diffusive magnetohydrodynamics.
For instance, the methods of \cite{FT} gave an opportunity to investigate
the Gevrey class regularity of the MHD solutions and to obtain the results \cite{Ki}
analogous to \cite{FT}.

The present paper has three goals:

\noindent
$\bullet$ to carry over the a priori bounds for arbitrary-order space derivatives
of solutions to the Navier--Stokes equation to space-periodic solutions
to the equations of diffusive magnetohydrodynamics;

\noindent
$\bullet$ to derive similar a priori bounds for arbitrary-order space derivatives
of the first-order time derivative of the Fourier--Galerkin approximants
and to prove that the bounds are admitted by weak solutions to the equations
of magnetohydrodynamics;

\noindent
$\bullet$ to reveal a link between these bounds and space analyticity of the MHD
solutions at almost all times.

They are achieved by following an original approach \cite{Gev} based
on a transformation of coefficients in the expansion of the solutions
in Fourier series in spatial variables. We introduce an auxiliary problem,
whose solutions are Fourier series involving the transformed coefficients;
an additional first-order pseudodifferential operator emerges in it. This
enables us to estimate a Gevrey class norm of the MHD solutions. The time-dependent
index of this norm, controlling the size in $\C$ (in the imaginary directions)
of the region of space analyticity upon complexification of the spatial
variables, is inversely proportional to a Sobolev norm of the solution to the
auxiliary system of equations. The estimate is global, i.e., applicable at all
times except for a set of Lebesgue measure zero, where the norm becomes
infinite. Finiteness of a Gevrey class norm of a solution implies that
its Fourier series converges as a geometric series, as well as the Fourier
series of its spatial derivatives. Following this observation, we construct
bounds for norms of arbitrary high spatial derivatives in terms of estimates of
a suitable norm of the solution and the common ratio of the geometric series.

The structure of the paper is as follows. In the next section we state
the problem, introduce the main equations to be investigated and set
the notation. In section~\ref{io} we follow \cite{FT} to show that space
analyticity sets in instantaneously, provided the initial data belongs
to the Sobolev space $H_s(\T)$. We are only interested in real analyticity.
We introduce in section~\ref{bo} the auxiliary system of equations and derive
an a priori bound of the energy type for its solutions. It is used
for construction of a priori bounds for Sobolev spaces and Wiener algebra norms
of weak solutions to the equations of magnetohydrodynamics in section~\ref{mhb},
and of the first-order time derivatives of the solutions in section~\ref{dt}.
While carrying a priori bounds for Fourier--Galerkin approximants over
to the weak solutions relies on standard arguments and is straightforward,
this is not the case of bounds for the time derivatives. They are justified
in section~\ref{bw}. We make the concluding remarks in the last section
of the paper. For the reader's convenience, our presentation is
reasonably detailed. The end of the proof of a lemma or theorem
is marked by the symbol $\RHD$.

This paper is dedicated to Professor Uriel Frisch on the occasion of his 80th
anniversary as a sign of appreciation of the Scientist and the Teacher.

\pagebreak
\section{Statement of the problem}\label{stp}

An electrically conducting fluid flow, whose velocity in the Eulerian
coordinates ${\bf x}\in\R$ is ${\bf V(x},t)$, in the presence of magnetic
field ${\bf B(x},t)$ satisfies the equations
\se{vb}\be{\partial{\bf V}\over\partial t}&=\nu\nabla^2{\bf V}
-({\bf V}\cdot\nabla){\bf V}+({\bf B}\cdot\nabla){\bf B}-\nabla P,\label{ve}\\
{\partial{\bf B}\over\partial t}&=\eta\nabla^2{\bf B}
+\nabla\times({\bf V}\times{\bf B}).\label{be}\end{align}
Here $P$ is the total pressure and $t$ is time. The first equation, \rf{ve}, is
the fluid momentum equation known as the Navier--Stokes equation, and the second
one, \rf{be}, is the magnetic induction equation.
We assume that the only external body force acting on the fluid is the magnetic
Lorentz force. (This assumption is made for the sake of simplicity only;
adding a prescribed space-analytic body force does not present any fundamental
mathematical difficulty, but makes the presentation more involved.) The flow
is supposed to be incompressible, and magnetic field is solenoidal:
\BE\nabla\cdot{\bf V}=\nabla\cdot{\bf B}=0.\EE{bvs}\end{subequations}
Initially (at $t=0$) the flow velocity ${\bf V}^{\rm(init)}$ and magnetic
field ${\bf B}^{\rm(init)}$ are prescribed.

We seek space-periodic solutions, the periodicity cell being a cube
$\T=[0,2\pi]^3$. Expanding the solution in Fourier series
\BE{\bf V}=\sum_{\bf n}\hV_{\bf n}\e^{\i\bf n\cdot x},\qquad
{\bf B}=\sum_{\bf n}\hB_{\bf n}\e^{\i\bf n\cdot x}\EE{Fs}
(where summation is over three-dimensional vectors $\bf n$ with integer
components), multiplying \rf{ve} and \rf{be} by $\e^{-\i\bf n\cdot x}$ and
integrating over $\T$ yields a system of ODEs for the Fourier coefficients
\se{wF}\be{\d\hV_{\bf n}\over\d t}&+\nu|{\bf n}|^2\hV_{\bf n}
=-\i\sum_{\bf k}\P\big((\hV_{\bf n-k}\!\cdot\!{\bf k})\hV_{\bf k}
-(\hB_{\bf n-k}\!\cdot\!{\bf k})\hB_{\bf k}\big),\label{wV}\\
{\d\hB_{\bf n}\over\d t}&+\eta|{\bf n}|^2\hB_{\bf n}
=\i\sum_{\bf k}{\bf n}\times(\hV_{\bf n-k}\times\hB_{\bf k}).
\label{wB}\end{align}\end{subequations}

Fourier--Galerkin approximants of solutions to \rf{vb} are truncated series
\BE{\bf V}^{(N)}=\sum_{|{\bf n}|\le N}\hV_{\bf n}^{(N)}\e^{\i\bf n\cdot x},\qquad
{\bf B}^{(N)}=\sum_{|{\bf n}|\le N}\hB_{\bf n}^{(N)}\e^{\i\bf n\cdot x}\EE{Fou}
(we set $\hV_{\bf n}^{(N)}=\hB_{\bf n}^{(N)}=0$ for $|{\bf n}|>N$). Fourier
coefficients of the approximants satisfy \rf{wF} for $|{\bf n}|\le N$.
Henceforth, we drop the superscript $(N)$ indicating the dependence
of the approximants on the resolution parameter $N$, but reinstate this
notation in section~\ref{bw}.

The fields $\bf V$ and $\bf B$ are assumed to be zero-mean,
\BE(2\pi)^{-3}\int_{\T}{\bf V}\,\d{\bf x}=(2\pi)^{-3}\int_{\T}{\bf B}\,\d{\bf x}=0
\quad\Leftrightarrow\quad\hV_0=\hB_0=0\EE{zm}
(note that in the course of temporal evolution due to equations \rf{vb} the
spatial means of $\bf V$ and $\bf B$ are conserved). They are real as long as
$$\overline{\hV_{\bf n}}=\hV_{\bf-n},\qquad\overline{\hB_{\bf n}}=\hB_{\bf-n}$$
(the bar denotes complex conjugation). The solenoidality conditions
\rf{bvs} reduce to the orthogonality
\BE\hV_{\bf n}\cdot{\bf n}=\hB_{\bf n}\cdot{\bf n}=0.\EE{so}
We denote by $\P$ the linear projection of a three-dimensional vector
on the plane normal to ${\bf n}\ne0$:
$$\P:{\bf f}\mapsto{\bf f}-{{\bf f}\cdot{\bf n}\over|{\bf n}|^2}{\bf n}.$$

Let $|\cdot|_p$ denote the norm in the functional Lebesgue space $L^p(\T)$,
$${|f|}_p=\Big((2\pi)^{-3}\int_{\T}|f|^p\,\d{\bf x}\Big)^{\!\!1/p},
\qquad p\ge1.$$
We denote by $\dot H_s(\T)$ the subspace of $H_s(\T)$ comprised
of space-periodic (with the periodicity cell $\T$) three-dimensional
zero-mean solenoidal vector fields equipped with the norm
$$\|{\bf w}\|_s=\Big(\sum_{\bf n}|{\bf n}|^{2s}|\widehat{\bf w}_{\bf n}|^2
\Big)^{\!\!1/2}\!\!\!,\quad\mbox{where~~}
{\bf w}=\sum_{{\bf n}\ne 0}\widehat{\bf w}_{\bf n}\e^{\i\bf n\cdot x}.$$
By the embedding theorem for Sobolev spaces (\cite{Pe}, see also, e.g.,
\cite{Li,AF,BL,Tr,Ma}),
for every positive $s<3/2$ there exists a constant $C_s$ such that each
function $f\in H_s(\T)$ of a three-dimensional space variable satisfies
the inequality
\BE{|f|}_{6/(3-2s)}\le C_s\|f\|_s.\EE{emb}

The Gevrey class norm is defined for $\sigma>0$ by the relation
$$\Gn{\bf w}\Gn^2_{\sigma,s}=\sum_{\bf n}|\widehat{\bf w}_{\bf n}|^2
\e^{2\sigma|\bf n|}|{\bf n}|^{2s}\quad\mbox{for~~}
{\bf w}=\sum_{\bf n}\widehat{\bf w}_{\bf n}\e^{\i\bf n\cdot x}.$$
If the norm of a field is finite, it is space-analytic, the size of the open
region of analyticity of the field in the imaginary directions for complex
$\bf x$ being at least $\sigma$. The inequality
\BE|{\bf n}|^a\e^{-b|{\bf n}|}\le(\e b/a)^{-a}\quad\mbox{for all}~~a>0,~~b>0,\EE{exp}
implies a relation between Gevrey norms of different indices:
\BE\Gn{\bf w}\Gn_{\sigma'\!,\,p}\le\Gn{\bf w}\Gn_{\sigma,s}(\e(\sigma-\sigma')/
(p-s))^{p-s}\quad\mbox{for}~~\sigma'<\sigma,~~s<p.\EE{gn}

\section{Instantaneous onset of space analyticity}\label{io}

In this section, we prove

{\it Theorem 1.} Let the initial data ${\bf V}^{\rm(init)}$ and
${\bf B}^{\rm(init)}$ at time $t=0$ belong to $\dot H_s(\T)$ for some $s>1/2$.
Then there exists $t_*>0$ such that the weak solution to the system of equations \rf{vb}
is space-analytic in the open interval $0<t<t_*$.

{\it Proof.} Following \cite{FT}, we set
\BE\hV_{\bf n}={\bf v_n}\e^{-\sigma|\bf n|t},\quad
\hB_{\bf n}={\bf b_n}\e^{-\sigma|{\bf n}|t},\EE{tsu}
where $\sigma<\min(\nu,\eta)$ is a strictly positive constant, and derive
a priori bounds on the interval $0<t<t_*$ for the modified solutions
$${\bf v}=\sum_{0\ne|{\bf n}|\le N}{\bf v_n}(t)\e^{\i\bf n\cdot x},\quad
{\bf b}=\sum_{0\ne|{\bf n}|\le N}{\bf b_n}(t)\e^{\i\bf n\cdot x}.$$
Substituting the truncated series \rf{Fou} and \rf{tsu} into \rf{wF} yields
\se{GF}\be{\d{\bf v_n}\over\d t}=&-(\nu|{\bf n}|^2-\sigma|{\bf n}|){\bf v_n}
-\i\sum_{\bf k}\e^{\sigma t(|{\bf n}|-|{\bf k}|-|{\bf n-k}|)}
\P\big(({\bf v_{n-k}}\!\cdot\!{\bf k}){\bf v_k}
-({\bf b_{n-k}}\!\cdot\!{\bf k}){\bf b_k}\big),\label{Ftv}\\
{\d{\bf b_n}\over\d t}=&-(\eta|{\bf n}|^2-\sigma|{\bf n}|){\bf b_n}
+\i{\bf n}\times\sum_{\bf k}\e^{\sigma t(|{\bf n}|-|{\bf k}|-|{\bf n-k}|)}
\,{\bf v_{n-k}}\times{\bf b_k}.\label{Ftb}\end{align}\end{subequations}

By the triangle inequality, the exponential
in the r.h.s.~of equations \rf{GF} does not exceed 1. For $1/2<s\le 1$,
we scalar multiply \rf{Ftv} and \rf{Ftb} by $|{\bf n}|^{2s}{\bf v_{-n}}$ and
$|{\bf n}|^{2s}{\bf b_{-n}}$, respectively, sum up the results over all
${\bf n}\ne 0$ (see \rf{zm}) and take into account the inequalities
\BE|{\bf v_{n-k}}\!\cdot\!{\bf k}|\le|{\bf k}|^\beta|{\bf n}|^{1-\beta}|{\bf v_{n-k}}|,\qquad
|{\bf b_{n-k}}\!\cdot\!{\bf k}|\le|{\bf k}|^\beta|{\bf n}|^{1-\beta}|{\bf b_{n-k}}|
\qquad\mbox{for any~}\beta,~~0\le\beta\le1\EE{ob}
(stemming from the orthogonality \rf{so}) and
\BE|{\bf n}|^s\le\max(1,2^{s-1})(|{\bf k}|^s+|{\bf n-k}|^s)
\quad\mbox{for any~}s\ge0.\EE{pin}
Choosing $\beta=s$, we find
\be&{1\over2}\,{\d\over\d t}(\|{\bf v}\|_s^2+\|{\bf b}\|_s^2)
+(\nu-\sigma)\|{\bf v}\|_{1+s}^2+(\eta-\sigma)\|{\bf b}\|_{1+s}^2\nonumber\\
&\le\sum_{\bf n,k}\Big(\big(|{\bf v_{n-k}}||{\bf v_k}|
+|{\bf b_{n-k}}||{\bf b_k}|\big)|{\bf k}|^s|{\bf v_{-n}}|
+\big(|{\bf k}|^s+|{\bf n-k}|^s\big)|{\bf v_{n-k}}||{\bf b_k}||{\bf b_{-n}}|
\Big)|{\bf n}|^{1+s}.\label{ene}\end{align}
In terms of scalar functions
\BE f^{\bf g}_q({\bf x},t)=\sum_{\bf n}|{\bf g_n}(t)||{\bf n}|^q
\e^{\i\bf n\cdot x}\quad\mbox{for an arbitrary~~}
{\bf g}({\bf x},t)=\sum_{\bf n}{\bf g_n}(t)\e^{\i\bf n\cdot x},\EE{scf}
the r.h.s.~of \rf{ene} can be expressed as
\ba&(2\pi)^{-3}\int_{\T}\Big(\big(f^{\bf v}_0f^{\bf v}_s
+f^{\bf b}_0f^{\bf b}_s\big)f^{\bf v}_{1+s}+\big(f^{\bf b}_0f^{\bf v}_s
+f^{\bf b}_sf^{\bf v}_0\big)f^{\bf b}_{1+s}\Big)\,\d{\bf x}\\
\shortintertext{and further bounded as follows:}
&\le\big(|f^{\bf v}_0|_{6/(3-2s)}|f^{\bf v}_s|_{3/s}
+|f^{\bf b}_0|_{6/(3-2s)}|f^{\bf b}_s|_{3/s}\big)|f^{\bf v}_{1+s}|_2
+\big(|f^{\bf b}_0|_{6/(3-2s)}|f^{\bf v}_s|_{3/s}
+|f^{\bf b}_s|_{3/s}|f^{\bf v}_0|_{6/(3-2s)}\big)|f^{\bf b}_{1+s}|_2\\
\shortintertext{(by H\"older's inequality)}
&\le C'_s\Big(\!\big(\|f^{\bf v}_0\|_s\|f^{\bf v}_s\|_{3/2-s}
\!+\|f^{\bf b}_0\|_s\|f^{\bf b}_s\|_{3/2-s}\big)\|f^{\bf v}_0\|_{1+s}
\!+\big(\|f^{\bf b}_0\|_s\|f^{\bf v}_s\|_{3/2-s}
\!+\|f^{\bf v}_0\|_s\|f^{\bf b}_s\|_{3/2-s}\big)\|f^{\bf b}_0\|_{1+s}\!\Big)\\
\shortintertext{(by the embedding theorem inequalities \rf{emb};
we have denoted $C'_s=C_sC_{3/2-s}$)}
&\le C'_s\Big(\|{\bf v}\|_s^{s+{1\over2}}\|{\bf v}\|_{1+s}^{{5\over2}-s}
\!+\!\|{\bf b}\|_s^{s+{1\over2}}\|{\bf b}\|_{1+s}^{{3\over2}-s}\|{\bf v}\|_{1+s}
\!+\!\|{\bf b}\|_s\|{\bf b}\|_{1+s}\|{\bf v}\|_s^{s-{1\over2}}\|{\bf v}\|_{1+s}^{{3\over2}-s}
\!+\!\|{\bf v}\|_s\|{\bf b}\|_s^{s-{1\over2}}\|{\bf b}\|_{1+s}^{{5\over2}-s}\Big)\\
\shortintertext{(since
$\|{\bf g}\|_{3/2}\le\|{\bf g}\|_s^{s-1/2}\|{\bf g}\|_{1+s}^{3/2-s}$
for any $\bf g$ by H\"older's inequality)}
&\le C'_s\Big(\!(5/2-s)\gamma(\|{\bf b}\|_{1+s}^2\!+\|{\bf v}\|_{1+s}^2)
+(s-1/2)\gamma^{-(5-2s)/(2s-1)}
(\|{\bf b}\|_s^{2(2s+1)/(2s-1)}\!+\|{\bf v}\|_s^{2(2s+1)/(2s-1)})\!\Big)\\
\shortintertext{(by Young's inequality; $\gamma>0$ is an arbitrary constant)}
&\le C'_s\Big(\!(5/2-s)\gamma(\|{\bf b}\|_{1+s}^2\!+\|{\bf v}\|_{1+s}^2)
+(s-1/2)\gamma^{-(5-2s)/(2s-1)}(\|{\bf b}\|_s^2\!+\|{\bf v}\|_s^2)^{(2s+1)/(2s-1)}\Big).
\end{align*}

Choosing now $\gamma>0$ such that $C'_s(5/2-s)\gamma\le\min(\nu,\eta)-\sigma$,
denoting $C''_s=C'_s\gamma^{-(5-2s)/(2s-1)}$, solving the inequality \rf{ene}
and applying \rf{tsu} yields
\BE\Gn{\bf V}\Gn_{\sigma t,s}^2+\Gn{\bf B}\Gn_{\sigma t,s}^2
=\|{\bf v}\|_s^2+\|{\bf b}\|_s^2\le((\|{\bf V}^{\rm(init)}\|_s^2
+\|{\bf B}^{\rm(init)}\|_s^2)^{-2/(2s-1)}-C''_st)^{-(s-1/2)}=q_s(t)\EE{hsb}
for
\BE t<t_*=(\|{\bf V}^{\rm(init)}\|_s^2+\|{\bf B}^{\rm(init)}
\|_s^2)^{-2/(2s-1)}/C''_s.\EE{tb}
(The initial data for the Fourier--Galerkin equations \rf{GF} should be used
in the r.h.s.~of \rf{hsb} and \rf{tb}, but we replace the norms
in \rf{hsb} and \rf{tb} by the norms of the initial data for the original
problem \rf{vb}, since the $H_s(\T)$ norms of the truncated initial conditions
monotonically increase with the resolution parameter $N$.) Thus,
the Fourier--Galerkin approximants \rf{Fou}
of solutions to \rf{vb} obey an a priori bound, that is independent of $N$.
Usual arguments show that they converge to a weak solution to the problem
\rf{vb} (see section \ref{bw}), that is space-analytic for $0<t<t_*$ even
if the initial condition is not, and hence on this time interval
it is strong and unique.~~$\RHD$

\section{An a priori bound for approximants of solutions to the auxiliary
problem}\label{bo}

We now consider the initial-value problem stated at $t=t_0$ such that
$0<t_0<t_*$. We have shown in the previous section that the ``initial'' fields
${\bf V}^{\rm(beg)}={\bf V}({\bf x},t_0)$ and
${\bf B}^{\rm(beg)}={\bf B}({\bf x},t_0)$ are space-analytic. By virtue of
\rf{hsb} and \rf{gn}, $\Gn{\bf V}^{\rm(beg)}\Gn_{\sigma'\!,3/2}^2
+\Gn{\bf B}^{\rm(beg)}\Gn_{\sigma'\!,3/2}^2\le q_{3/2}(t_0)<\infty$
for any $\widetilde\sigma<\sigma t_0$. Constructions of the present section
are based on this property of the initial data and otherwise do not rely
on the results of the previous section: It suffices to assume that the norms
$\Gn\!\cdot\!\Gn_{\sigma'\!,3/2}$ of the data at $t=t_0$ are finite for some
$\widetilde\sigma>0$ and uniformly in $N$ bounded for the finite-space
Fourier--Galerkin approximants, and consider the initial-value
problem paying no attention to the prior existence of the solution
to \rf{vb} for $0\le t\le t_0$. We assume
$\|{\bf V}^{\rm(beg)}\|_{3/2}+\|{\bf B}^{\rm(beg)}\|_{3/2}>0$.

\subsection{A transformation of solutions to \rf{vb} and the auxiliary
system of equations}\label{tra}

Following \cite{Gev}, we transform the Fourier coefficients
\BE\hV_{\bf n}=\wv_{\bf n}\e^{-\delta\Phi|{\bf n}|},\quad
\hB_{\bf n}=\wb_{\bf n}\e^{-\delta\Phi|{\bf n}|}\EE{svb}
of the truncated Fourier--Galerkin approximants \rf{Fou} of a solution
to \rf{vb}; here we have denoted
\BE\Phi=(1+\|\wv\|_{3/2}^2+\|\wb\|_{3/2}^2)^{-1/2},\quad
\wv=\sum_{0\ne|{\bf n}|\le N}\wv_{\bf n}(t)\e^{\i\bf n\cdot x},\quad
\wb=\sum_{0\ne|{\bf n}|\le N}\wb_{\bf n}(t)\e^{\i\bf n\cdot x}\EE{wvb}
and $\delta>0$ is a constant. Substituting
the Fourier series \rf{Fou} and \rf{svb} into \rf{wF} yields
\se{GFm}\be{\d\wv_{\bf n}\over\d t}&+\nu|{\bf n}|^2\wv_{\bf n}
-\delta|{\bf n}|\wv_{\bf n}\,{\d\Phi\over\d t}
=-\i\sum_{\bf k}\e^{\delta\Phi(|{\bf n}|-|{\bf k}|-|{\bf n-k}|)}
\P\big((\wv_{\bf n-k}\!\cdot\!{\bf k})\wv_{\bf k}
-(\wb_{\bf n-k}\!\cdot\!{\bf k})\wb_{\bf k}\big),\label{mv}\\
{\d\wb_{\bf n}\over\d t}&+\eta|{\bf n}|^2\wb_{\bf n}
-\delta|{\bf n}|\wb_{\bf n}\,{\d\Phi\over\d t}
=\i\sum_{\bf k}\e^{\delta\Phi(|{\bf n}|-|{\bf k}|-|{\bf n-k}|)}\,
{\bf n}\times(\wv_{\bf n-k}\times\wb_{\bf k}).
\label{mb}\end{align}\end{subequations}

The system of ODEs \rf{GFm} is satisfied
by the Fourier--Galerkin approximants of solutions $\wv,\wb$ to the system
of pseudodifferential equations which we call an {\it auxiliary problem}:
\se{tr}\be{\partial\wv\over\partial t}-\D\wv\,{\d\over\d t}\ln\Phi
&=\nu\nabla^2\wv+\e^{\D}\Big(
-\big(\big(\e^{-\D}\wv\big)\cdot\nabla\big)\big(\e^{-\D}\wv\big)
+\big(\big(\e^{-\D}\wb\big)\cdot\nabla\big)\big(\e^{-\D}\wb\big)
\Big)-\nabla\widetilde p,\label{trv}\\
{\partial\wb\over\partial t}-\D\wb\,{\d\over\d t}\ln\Phi
&=\eta\nabla^2\wb+\nabla\times\Big(\e^{\D}\Big(\big(\e^{-\D}\wv\big)
\times\big(\e^{-\D}\wb\big)\Big)\Big),\label{trb}\\
\nabla\cdot\wv&=\nabla\cdot\wb=0\label{trs}\end{align}\end{subequations}
where $\D=\delta\Phi(-\nabla^2)^{1/2}$ is defined in the subspace of zero-mean
vector fields of the Lebesgue space $L^2(\T)$. While not very illuminating
in the present setup, the equations in this form may be useful when considering
the problem with appropriate boundary conditions in a finite fluid domain.
For $\delta=0$, they reduce to the original equations of magnetohydrodynamics
\rf{vb}.

To render \rf{GFm} as an explicit system of ODEs,
we scalar multiply \rf{mv} and \rf{mb} by $|{\bf n}|^3{\wv_{\bf-n}}$
and $|{\bf n}|^3{\wb_{\bf-n}}$, respectively, sum up the results over all
${\bf n}\ne 0$ (see \rf{zm}) and obtain
\BE{1\over2}\,(1+\delta\Phi^3(\|\wv\|_2^2+\|\wb\|_2^2))
{\d\over\d t}(\|\wv\|_{3/2}^2+\|\wb\|_{3/2}^2)
+\nu\|\wv\|_{5/2}^2+\eta\|\wb\|_{5/2}^2=\i\Sigma_3,\EE{der}
where we have denoted
\be\Sigma_p&=\!\sum_{\bf n,k}|{\bf n}|^p\e^{\delta\Phi(|{\bf n}|-|{\bf k}|-|{\bf n-k}|)}
\Big(\!-(\wv_{\bf n-k}\!\cdot\!{\bf k})\wv_{\bf k}\!\cdot\!\wv_{\bf-n}
+(\wb_{\bf n-k}\!\cdot\!{\bf k})\wb_{\bf k}\!\cdot\!\wv_{\bf-n}
+({\bf n}\times(\wv_{\bf n-k}\times\wb_{\bf k}))\!\cdot\!\wb_{\bf-n}\Big)\nonumber\\
&=\!\sum_{\bf n,k}|{\bf n}|^p\e^{\delta\Phi(|{\bf n}|-|{\bf k}|-|{\bf n-k}|)}
\Big(-(\wv_{\bf n-k}\!\cdot\!{\bf k})(\wv_{\bf k}\!\cdot\!\wv_{\bf-n}
\!+\wb_{\bf k}\!\cdot\!\wb_{\bf-n})
+(\wb_{\bf n-k}\!\cdot\!{\bf k})(\wb_{\bf k}\!\cdot\!\wv_{\bf-n}
\!+\wv_{\bf k}\!\cdot\!\wb_{\bf-n})\Big)\label{sp}\end{align}
(here and in what follows, we use the orthogonality relations
$\wv_{\bf k}\cdot{\bf k}=\wb_{\bf k}\cdot{\bf k}=0$ stemming from \rf{so},
and swap the indices of summation $\bf k$ and $\bf n-k$ in some terms
when it is convenient to rearrange the sums).
Finding $\d\Phi/\d t=-(\Phi^3/2){\d\over\d t}(\|\wv\|_{3/2}^2\!+\|\wb\|_{3/2}^2)$
from \rf{der} and substituting into \rf{GFm} yields the desired explicit
system of ODEs, for which we now need to supply the initial conditions.

The transformed solutions can be constructed for both the truncated sums \rf{Fou}
and infinite series~\rf{Fs}. By virtue of \rf{svb}, the harmonics
$\wv^{\rm(beg)}_{\bf n}=\hV^{\rm(beg)}_{\bf n}\e^{\delta\Phi|{\bf n}|},~~
\wb^{\rm(beg)}_{\bf n}=\hB^{\rm(beg)}_{\bf n}\e^{\delta\Phi|{\bf n}|}$
are available at $t=t_0$, if the value of the parameter $\Phi$ is known.
These relations imply
\BE\Theta(\Phi(t_0))=0,\quad\mbox{where}\quad\Theta(\Phi)
=\Gn{\bf V}^{\rm(beg)}\Gn_{\delta\Phi,3/2}^2
+\Gn{\bf B}^{\rm(beg)}\Gn_{\delta\Phi,3/2}^2+1-\Phi^{-2}.\EE{eqn}
We regard \rf{eqn} as an equation in $\Phi(t_0)$.
Let $R$ denote the size of the region of analyticity, defined as the infimum
of such $r$ that $\Gn{\bf V}^{\rm(beg)}\Gn_{r,3/2}^2
+\Gn{\bf B}^{\rm(beg)}\Gn_{r,3/2}^2=\infty$.
By the results of the previous section, $R\ge\sigma t_0>0$; evidently,
$R=\infty$ for truncated series \rf{Fs}. We note that $0<\Phi\le1$ and
$\Theta(\Phi)$ increases monotonically from $\Theta(0)=-\infty$ to
$\Theta(1)=\Gn{\bf V}^{\rm(beg)}\Gn_{\delta,3/2}^2
+\Gn{\bf B}^{\rm(beg)}\Gn_{\delta,3/2}^2>0$
if $R>\delta$, or to $\Theta(R/\delta)=\infty$ otherwise.
Hence, in both cases \rf{eqn} has a unique solution $\Phi(t_0)$.
Now the initial data $\wv^{\rm(beg)}_{\bf n}$ and $\wb^{\rm(beg)}_{\bf n}$
at $t=t_0$ are fully determined.

\subsection{The ``energy'' bound for the transformed solutions}\label{eb}

Thus, for a given resolution parameter $N$, $\wv_{\bf n}$ and $\wb_{\bf n}$ can
be found for any $t>t_0$ as a solution to an explicit finite system of ODEs,
provided
it does not blow up at a finite time. Theorem 2 rules out this possibility.

{\it Theorem 2.} Suppose
\BE\delta\le(18\sqrt2C'_{1/2})^{-1}\min(\nu,\eta),\EE{del}
where $C'_{1/2}=C_{1/2}C_1$ (see \rf{emb}). Solutions $\wv,\wb$
to the auxiliary problem \rf{wvb} obey an a priori bound
\BE{9\over4}\left.(\|\wv\|_0^2+\|\wb\|_0^2)\right|_{t=T}
+\int_{t_0}^T\!\!\left(\nu\big(\|\wv\|_1^2+4\delta^2\Phi^2
\|\wv\|_2^2\big)\!+\eta\big(\|\wb\|_1^2+4\delta^2\Phi^2
\|\wb\|_2^2\big)\!\right)\!\d t\le9Q\EE{bou}
for all $T\ge t_0$, where
$$Q=\Big({1\over2}(\|\wv\|_0^2+\|\wb\|_0^2)
-\delta\Phi(\|\wv\|_{1/2}^2+\|\wb\|_{1/2}^2)
+\delta^2\Phi^2(\|\wv\|_1^2+\|\wb\|_1^2)
+{2\delta^3\over3}(\Phi^3-3\Phi+2)\!\Big)\Big|_{t=t_0}$$
depends only on the initial conditions $\hV^{\rm(beg)},\hB^{\rm(beg)}$
at $t=t_0$ and the parameter $\delta$.

{\it Proof.} We consider an analogue of the energy equation for \rf{GFm}.
Because of the presence of new terms in \rf{GFm}, we scalar multiply \rf{mv}
and \rf{mb} by $\wv_{\bf-n}(1-2\delta|{\bf n}|\Phi+2\delta^2|{\bf n}|^2\Phi^2)$
and $\wb_{\bf-n} (1-2\delta|{\bf n}|\Phi+2\delta^2|{\bf n}|^2\Phi^2)$,
respectively. Summing up the results over $\bf n$ yields
\be&{\d\over\d t}\left({1\over2}(\|\wv\|_0^2+\|\wb\|_0^2)
-\delta\Phi(\|\wv\|_{1/2}^2+\|\wb\|_{1/2}^2)
+\delta^2\Phi^2(\|\wv\|_1^2+\|\wb\|_1^2)
+2\delta^3\left({\Phi^3\over3}-\Phi\right)\!\!\right)\nonumber\\
&+\nu\big(\|\wv\|_1^2-2\delta\Phi\|\wv\|_{3/2}^2
+2\delta^2\Phi^2\|\wv\|_2^2\big)+\eta\big(\|\wb\|_1^2
-2\delta\Phi\|\wb\|_{3/2}^2+2\delta^2\Phi^2\|\wb\|_2^2\big)\nonumber\\
&=\i\left(\Sigma_0-2\delta\Phi\Sigma_1+2\delta^2\Phi^2\Sigma_2
\right)\label{lhs}\end{align}
(see \rf{sp}). The terms constituting the l.h.s.~of \rf{lhs} are bounded as follows:
$${1\over2}(\|\wv\|_0^2+\|\wb\|_0^2)-\delta\Phi(\|\wv\|_{1/2}^2+\|\wb\|_{1/2}^2)
+\delta^2\Phi^2(\|\wv\|_1^2+\|\wb\|_1^2)\ge{1\over4}(\|\wv\|_0^2+\|\wb\|_0^2),$$
since $\|\wv\|_{1/2}^2+\|\wb\|_{1/2}^2\le\|\wv\|_0\|\wv\|_1+\|\wb\|_0\|\wb\|_1$;
$$0>2\delta^3(\Phi^3/3-\Phi)\ge-4\delta^3/3,$$
since on the interval $0\le\Phi\le1$ the function $\Phi^3/3-\Phi$ monotonically
decreases;
$$\|\wv\|_1^2-2\delta\Phi\|\wv\|_{3/2}^2+2\delta^2\Phi^2\|\wv\|_2^2\ge
{1\over3}\|\wv\|_1^2+{\delta^2\Phi^2\over2}\|\wv\|_2^2,$$
since $\|\wv\|_{3/2}^2\le\|\wv\|_1\|\wv\|_2$ and due to the elementary
inequality $2|ab|\le\kappa^{-1}a^2+\kappa b^2$ that holds true for any $a,b$
and $\kappa>0$; similarly,
$$\|\wb\|_1^2-2\delta\Phi\|\wb\|_{3/2}^2+2\delta^2\Phi^2\|\wb\|_2^2\ge
{1\over3}\|\wb\|_1^2+{\delta^2\Phi^2\over2}\|\wb\|_2^2.$$

The sums $\Sigma_p$ in the r.h.s.~of \rf{lhs} are bounded by essentially
different procedures depending on whether $p$ vanishes or it is positive.
We note that the exponential does not exceed 1, since the exponent
in the r.h.s.~is negative. For $p=1,2$, \rf{sp} implies
\be|\Sigma_p|&\le\sum_{\bf n,k}|{\bf n}|^p|{\bf k}|
\left(|\wv_{\bf-n}|\big(|\wv_{\bf k}||\wv_{\bf n-k}|
+|\wb_{\bf k}||\wb_{\bf n-k}|\big)+|\wb_{\bf-n}|\big(
|\wb_{\bf n-k}||\wv_{\bf k}|+|\wv_{\bf n-k}||\wb_{\bf k}|
\big)\!\right)\label{s1}\\
&=(2\pi)^{-3}\int_{\T}\Big(f^{\wv}_p\big(f^{\wv}_1f^{\wv}_0
+f^{\wb}_0f^{\wb}_1\big)+f^{\wb}_p\big(f^{\wb}_0f^{\wv}_1
+f^{\wb}_1f^{\wv}_0\big)\Big)\,\d{\bf x}\nonumber\\
\shortintertext{(see \rf{scf})}
&\le|f^{\wv}_p|_2\big(|f^{\wv}_1|_3|f^{\wv}_0|_6
+|f^{\wb}_1|_3|f^{\wb}_0|_6\big)+|f^{\wb}_p|_2\big(
|f^{\wv}_1|_3|f^{\wb}_0|_6+|f^{\wb}_1|_3|f^{\wv}_0|_6\big)\nonumber\\
\shortintertext{(by H\"older's inequality)}
&\le C'_{1/2}\Big(\|\wv\|_p\big(\|\wv\|_{3/2}\|\wv\|_1
+\|\wb\|_{3/2}\|\wb\|_1\big)+\|\wb\|_p\big(
\|\wv\|_{3/2}\|\wb\|_1+\|\wb\|_{3/2}\|\wv\|_1\big)\!\Big)\nonumber\\
\shortintertext{(by the embedding theorem inequalities \rf{emb})}
&\le C'_{1/2}\sqrt2\big(\|\wv\|_p^2+\|\wb\|_p^2\big)^{1/2}
\big(\|\wv\|_{3/2}^2+\|\wb\|_{3/2}^2\big)^{1/2}
\big(\|\wv\|_1^2+\|\wb\|_1^2\big)^{1/2}\nonumber\end{align}
(by the Cauchy--Schwarz inequality). For $p=2$, this implies
$$2\delta^2\Phi^2|\Sigma_2|\le\delta C'_{1/2}\sqrt2\left(\|\wv\|_1^2
+\|\wb\|_1^2+\delta^2\Phi^2(\|\wv\|_2^2+\|\wb\|_2^2)\right).$$

For $p=0$, we symmetrise \rf{sp} by changing the indices of summation
$\bf n\to-k$ and $\bf k\to-n$ and using the solenoidality conditions~\rf{so}:
\ba\Sigma_0&={1\over2}\sum_{\bf n,k}\left(\e^{\delta\Phi(|{\bf n}|-|{\bf k}|-|{\bf n-k}|)}
-\e^{\delta\Phi(|{\bf k}|-|{\bf n}|-|{\bf n-k}|)}\right)\\
&\phantom{=}~\times\left(-(\wv_{\bf n-k}\!\cdot\!{\bf k})
(\wv_{\bf k}\!\cdot\!\wv_{\bf-n}+\wb_{\bf k}\!\cdot\!\wb_{\bf-n})
+(\wb_{\bf n-k}\!\cdot\!{\bf k})(\wb_{\bf k}\!\cdot\!\wv_{\bf-n}
+\wv_{\bf k}\!\cdot\!\wb_{\bf-n})\right).\\
\shortintertext{Evidently, $|\e^{\kappa_1}-\e^{\kappa_2}|\le|\kappa_1-\kappa_2|$
for all $\kappa_1<0$ and $\kappa_2<0$, implying}
|\Sigma_0|&\le\delta\Phi\sum_{\bf n,k}|{\bf n-k}||{\bf k}|
\left(|\wv_{\bf n-k}|(|\wv_{\bf k}||\wv_{\bf-n}|+|\wb_{\bf k}||\wb_{\bf-n}|)
+|\wb_{\bf n-k}|(|\wv_{\bf k}||\wb_{\bf-n}|+|\wb_{\bf k}||\wv_{\bf-n}|)
\right)\\
&\le\delta C'_{1/2}\sqrt2(\|\wv\|_1^2+\|\wb\|_1^2),\end{align*}
because the sum in the middle is identical to \rf{s1} for $p=1$.

Integrating \rf{lhs} in time and applying the above inequalities yields
\ba&{1\over4}\left.(\|\wv\|_0^2+\|\wb\|_0^2)\right|_{t=T}
\!+\!\int_{t_0}^T\!\!\left(
\nu\left({1\over3}\|\wv\|_1^2+{\delta^2\Phi^2\over2}\|\wv\|_2^2\right)
\!+\eta\left({1\over3}\|\wb\|_1^2+{\delta^2\Phi^2\over2}\|\wb\|_2^2\right)
\!\!\right)\!\d t\\
&\le Q+\delta C'_{1/2}\sqrt2\int_{t_0}^T\left(
4(\|{\bf v}\|_1^2+\|{\bf b}\|_1^2)+\delta^2\Phi^2(\|{\bf v}\|_2^2
+\|{\bf b}\|_2^2)\right)\!\d t.\end{align*}
Applying now the condition \rf{del}, we obtain the inequality \rf{bou}
as required.~~$\RHD$

\pagebreak
Compared to the usual energy bound, we have thus obtained a new bound
$$\int_{t_0}^T\!\Phi^2(\|\wv\|_2^2+\|\wb\|_2^2)\,\d t\le
9Q\left/(4\delta^2\min(\nu,\eta))\right.$$
for the Fourier--Galerkin approximants of solutions to the auxiliary problem
\rf{tr}. Although it is uniform over the resolution parameter $N$,
a further effort is required to deduce from it a bound
for weak solutions to the equations \rf{vb}. This ``bonus'' bound is due to the
presence of the first-order dissipative operator in the modified equation, that
emerges upon the transformation \rf{tsu}. An~operator of this type
was originally employed in the study of the ``lake'' equation in \cite{LO},
where a time dependence of the index of the respective Gevrey class norm
of the solution was assumed; our transformation \rf{tsu} also
introduces such a dependence, but a different one.

\section{A priori bounds for approximants of solutions to the system \rf{vb}}\label{mhb}

Here we use Theorem 2 for constructing bounds for the Fourier--Galerkin
approximants of solutions to the MHD system of equations \rf{vb} in Sobolev and
Wiener algebra norms, that are uniform in the truncation parameter $N$.
They feature the same exponents $\alpha_s$
as those considered in \cite{FGT}.\break A similar approach was entertained
in \cite{OT}, where bounds for algebraic decay of high-order derivatives
of strong solutions to the unforced Navier--Stokes equations in $R^n$ were
constructed by bounding a single Gevrey class seminorm of the solutions.

{\it Theorem 3.} For $\alpha_s=2/(2s-1),~\gamma_s=2/s$ and any $T>t_0$,
\se{zz}\be\int_{t_0}^T\!\!\big(\|{\bf V}\|_s^2+\|{\bf B}\|_s^2\big)^{\alpha_s/2}\d t
&\le\widetilde Q_s\mbox{~~for~~}s>1,\label{qq}\\
\int_{t_0}^T\!\!\big(\|{\bf V}\|_s^2+\|{\bf B}\|_s^2\big)^{\gamma_s/2}\d t
&\le\widetilde Q_s\mbox{~~for~~}0<s\le1,\label{sa}\\
\int_{t_0}^T\!\Big(\max_{\T}\big(|(-\nabla^2)^{s/2}{\bf V}|
+|(-\nabla^2)^{s/2}{\bf B}|\big)\!\Big)^{\alpha_{s+3/2}}\d t
&\le\widetilde Q^W_s\mbox{~~for~~}s>-1/2,\label{was}
\end{align}\end{subequations}
where $\widetilde Q_s$ and $\widetilde Q^W_s$ depend only on the initial
conditions $\hV^{\rm(beg)},\hB^{\rm(beg)}$ at $t=t_0$ and parameters
$\delta,s,\nu$ and~$\eta$. For $0<s\le1$, $\widetilde Q_s$ are
independent of $T\ge t_0$, while $\widetilde Q_s$ for $s>1$ and
$\widetilde Q^W_s$ are sublinear functions of $T-t_0$.

{\it Proof} differs in details when the inequalities for Sobolev and Wiener
algebra norms are considered. It is presented in the next two sections.

\subsection{Bounds in the Sobolev space norms}\label{bSo}

By the Cauchy--Schwarz inequality,
$$\|\wv\|_{3/2}^2+\|\wb\|_{3/2}^2\le\|\wv\|_1\|\wv\|_2+\|\wb\|_1\|\wb\|_2
\le(\|\wv\|_1^2+\|\wb\|_1^2)^{1/2}(\|\wv\|_2^2+\|\wb\|_2^2)^{1/2}.$$
This implies
$$(\|\wv\|_{3/2}^2+\|\wb\|_{3/2}^2)^{1/2}\le\left\{\begin{array}{ll}
\sqrt2\Phi(\|\wv\|_1^2+\|\wb\|_1^2)^{1/2}(\|\wv\|_2^2+\|\wb\|_2^2)^{1/2},
~&\mbox{if~}\|\wv\|_{3/2}^2+\|\wb\|_{3/2}^2\ge1,\\
\sqrt2\Phi(\|\wv\|_2^2+\|\wb\|_2^2)^{1/2},
~&\mbox{if~}\|\wv\|_{3/2}^2+\|\wb\|_{3/2}^2\le1.\end{array}\right.$$
Integrating in time and using \rf{bou} yields
\be\int_{t_0}^T\!\!(\|\wv\|_{3/2}^2&+\|\wb\|_{3/2}^2)^{1/2}\,\d t
\le\left(2\int_{t_0}^T\!\!\Phi^2(\|\wv\|_2^2+\|\wb\|_2^2)\,\d
t~\max\!\left(\int_{t_0}^T\!\!(\|\wv\|_1^2+\|\wb\|_1^2)\,\d t,T-t_0\!\right)
\!\!\right)^{\!\!\!1/2}\nonumber\\
&\le Q'=\delta^{-1}\max\Big(9Q\big/(\sqrt2\,\min(\nu,\eta)),
\big(9Q(T-t_0)\big/(2\,\min(\nu,\eta))\big)^{1/2}\Big).\label{Qp}\end{align}
Hence, for $s\ge1$ and $\alpha_s=2/(2s-1)$, the inequality
\BE(a+b)^c\le a^c+b^c\quad\mbox{for all}~~a>0,~~b>0,~~0\le c\le1,\EE{ab}
H\"older's inequality, \rf{Qp} and \rf{bou} imply
\be\int_{t_0}^T\!\!\big(&\|{\wv}\|_1^2+\|{\wb}\|_1^2\big)^{\alpha_s/2}\,
\Phi^{-(s-1)\alpha_s}\d t\le\int_{t_0}^T\!\!\big(\|{\wv}\|_1^2
+\|{\wb}\|_1^2\big)^{\alpha_s/2}\Big(1+\big(\|{\wv}\|_{3/2}^2
+\|{\wb}\|_{3/2}^2\big)^{(s-1)\alpha_s/2}\Big)\d t\nonumber\\
&\le Q''_s=\!\Big({9Q\over\min(\nu,\eta)}\Big)^{\!\alpha_s/2}
\big((T-t_0)^{1-\alpha_s/2}+(Q')^{1-\alpha_s/2}\big).\label{iti}\end{align}

We can now establish a priori bounds for the truncated Fourier--Galerkin
approximants \rf{Fou} of solutions to the original equations
of magnetohydrodynamics \rf{vb}. By \rf{svb} and the inequality \rf{exp},
$$\|{\bf V}\|_s^2+\|{\bf B}\|_s^2=\sum_{{\bf n}\ne0}|{\bf n}|^{2s}
\e^{-2\delta\Phi|{\bf n}|}\big(|\wv_{\bf n}|^2+|\wb_{\bf n}|^2\big)
\le\sum_{{\bf n}\ne0}|{\bf n}|^2\big(|\wv_{\bf n}|^2+|\wb_{\bf n}|^2\big)
\left({s-1\over\e\delta\Phi}\right)^{\!\!2(s-1)},$$
and hence applying \rf{iti} proves \rf{qq}:
$$\int_{t_0}^T\!\!\big(\|{\bf V}\|_s^2+\|{\bf B}\|_s^2\big)^{\alpha_s/2}\d t
\le\left(\!{s-1\over\e\delta}\right)^{\!\!(s-1)\alpha_s}Q''_s
\mbox{~~for~~}s>1.$$

For $0<s\le1$, the exponents are obtained by interpolating between
the endpoints of the interval, and thus constitute a different family.
Young's inequality and the energy inequality for \rf{vb} prove~\rf{sa}:
\ba\int_0^T\!\!\big(\|{\bf V}\|_s^2+\|{\bf B}\|_s^2\big)^{1/s}\,\d t
\le&\int_0^T\!\!(\|{\bf V}\|_1^2+\|{\bf B}\|_1^2)\,\d t\,
\max_{0\le t\le T}\big(\|{\bf V}\|_0^2+\|{\bf B}\|_0^2\big)^{(1-s)/s}\\
\le&\,(2\min(\nu,\eta))^{-1}(\|{\bf V}^{\rm(init)}\|_0^2
+\|{\bf B}^{\rm(init)}\|_0^2)^{1/s}.\end{align*}

{\it Corollary}. For any $p\ge2$, $s\ge3/p-1/2$ and $T>t_0$,
$$\int_{t_0}^T\!\!\big(\|{\bf V}\|_{p,s}
+\|{\bf B}\|_{p,s}\big)^{\alpha_{p,s}}\d t<\infty,$$
where we have denoted $\alpha_{p,s}=p/(p(s+1)-3)$ and
$\|{\bf f}\|_{p,s}=|(-\nabla^2)^{s/2}{\bf f}|_p$.

{\it Proof.} By \rf{emb}, ${|\bf f|}_p\!\le C_{3/2-3/p}\|{\bf f}\|_{3/2-3/p}$.
We apply this inequality to $(-\nabla^2)^{s/2}{\bf V}$ and
$(-\nabla^2)^{s/2}{\bf B}$, and note $\alpha_{p,s}=\alpha_{s+3/2-3/p}$.
By \rf{qq},
$$\int_{t_0}^T\!\!\big(\|{\bf V}\|_{p,s}
+\|{\bf B}\|_{p,s}\big)^{\alpha_{p,s}}\d t\le C_{3/2-3/p}^{\,\alpha_{p,s}}
\int_{t_0}^T\!\!\big(\|{\bf V}\|_{s+3/2-3/p}
+\|{\bf B}\|_{s+3/2-3/p}\big)^{\alpha_{s+3/2-3/p}}\d t<\infty.\quad\rhd$$
For solutions to the Navier--Stokes equation, this was proven in \cite{G4}.
Using \rf{sa}, it is easy to derive the analogous exponents for the case
$0<s+3/2-3/p\le1$.

\subsection{A priori bounds for the Wiener algebra norm}\label{bWi}

We finish here the proof of Theorem 3.

The bound \rf{was} follows from a bound for the Wiener algebra norm
of the fields $(-\nabla^2)^{s/2}{\bf V}$ and $(-\nabla^2)^{s/2}{\bf B}$.
The norm of a field is defined as the sum of absolute values of its Fourier
coefficients, i.e., a field has a finite Wiener algebra norm whenever
its Fourier series converges absolutely; obviously, the norm bounds
the field's maximum. (Applying this Banach space proved useful,
for instance, for estimating the dissipation length scale for turbulence
\cite{bis} and for showing time analyticity of solutions to the Euler equation
in Lagrangian coordinates \cite{uf1,uf2}.) The proof exploits the

{\it Lemma.} For any $\Phi,a$ and $p$ such that $0<\Phi\le1$, $a>0$ and $p>-3$,
$$\sum_{{\bf n}\ne0}\e^{-a\Phi|\bf n|}|{\bf n}|^p\le C_{p,a}^2\,\Phi^{-(p+3)},$$
where constants $C_{p,a}$ depend on $p$ and $a$, but not on $\Phi$.

{\it Proof.} Let $K_{\bf n}$ denote the cube
$\big\{{\bf k}\,\big||n_i-k_i|\le1/2\big\}$. Then
\ba\sum_{{\bf n}\ne0}\e^{-a\Phi|\bf n|}|{\bf n}|^p&=
\sum_{{\bf n}\ne0}\,\int_{K_{\bf n}}\e^{-a\Phi|\bf n|}|{\bf n}|^p\,\d{\bf k}\\
&\le\sum_{{\bf n}\ne0}\,\sup_{{\bf k}\in K_{\bf n}}\e^{a\Phi(|\bf k|-|\bf n|)}\,
\sup_{{\bf k}\in K_{\bf n}}(|{\bf n}|/|{\bf k}|)^p
\int_{K_{\bf n}}\e^{-a\Phi|\bf k|}|{\bf k}|^p\,\d{\bf k}\\
&\le\e^{a\Phi\sqrt3/2}\,2^{|p|}\int_{\R}\e^{-a\Phi|\bf k|}|{\bf k}|^p\,\d{\bf k}\\
&\le\e^{a\sqrt3/2}\,2^{|p|}\,4\pi\int_0^\infty\!\e^{-a\Phi\rho}\rho^{\,p+2}\,\d\rho\\
&=C_{p,a}^2\,\Phi^{-(p+3)},\end{align*}
where
$$C_{p,a}=\Big(4\pi\e^{a\sqrt3/2}\,2^{|p|}\,a^{-(p+3)}
\int_0^\infty\!\e^{-\rho}\rho^{\,p+2}\,\d\rho\Big)^{\!\!1/2}.~~\rhd$$
By \rf{svb}, the Cauchy--Schwarz inequality and Lemma, for $s>-1/2$,
\ba|(-\nabla^2)^{s/2}{\bf V}|+|(-\nabla^2)^{s/2}{\bf B}|
&=\Big|\sum_{{\bf n}\ne0}|{\bf n}|^s
\wv_{\bf n}\e^{-\delta\Phi|{\bf n}|+\i\bf n\cdot x}\Big|
+\Big|\sum_{{\bf n}\ne0}|{\bf n}|^s
\wb_{\bf n}\e^{-\delta\Phi|{\bf n}|+\i\bf n\cdot x}\Big|\\
&\le\sum_{{\bf n}\ne0}|{\bf n}|^s|\wv_{\bf n}|\e^{-\delta\Phi|{\bf n}|}
+\sum_{{\bf n}\ne0}|{\bf n}|^s|\wb_{\bf n}|\e^{-\delta\Phi|{\bf n}|}\\
&\le\Big(\!\Big(\sum_{{\bf n}\ne0}|{\bf n}|^2|\wv_{\bf n}|^2\Big)^{\!\!1/2}
\!+\Big(\sum_{{\bf n}\ne0}|{\bf n}|^2|\wb_{\bf n}|^2\Big)^{\!\!1/2}\,\Big)\!
\Big(\sum_{{\bf n}\ne0}|{\bf n}|^{2s-2}\e^{-2\delta\Phi|{\bf n}|}\Big)^{\!\!1/2}\\
&\le\sqrt2(\|\wv\|_1^2+\|\wb\|_1^2)^{1/2}\,C_{2s-2,2\delta}\,\Phi^{-(s+1/2)}.\end{align*}
Thus, application of \rf{iti} upon changing $s\to s+3/2$ for $s>-1/2$
establishes \rf{was}:
$$\int_{t_0}^T\!\!\Big(\max_{\T}\big(|(-\nabla^2)^{s/2}{\bf V}|
+|(-\nabla^2)^{s/2}{\bf B}|\big)\!\Big)^{\alpha_{s+3/2}}\,\d t
\le(\sqrt2C_{2s-2,2\delta})^{\alpha_{s+3/2}}\,Q''_{s+3/2}.~~\RHD$$

The bound \rf{was} was proven for solutions of the Navier--Stokes equation
for $s=0$ in \cite{FGT} (the authors attribute the proof to L.~Tartar
\cite{Ta}), and for $s>0$ in \cite{G4}.

\section{A priori bounds for time derivatives of solutions to the system \rf{vb}}\label{dt}

Similar bounds for higher-order norms of $\partial{\bf V}\!/\partial t$ and
$\partial{\bf B}/\partial t$ can now be constructed by using space analyticity
of the solutions to \rf{vb}. An alternative derivation based on bounds \rf{qq}
for the solutions is presented in section~\ref{alt}.

{\it Theorem 4}. Time derivatives of the solutions to the system of equations of
magnetohydrodynamics \rf{vb} satisfy the a priori inequalities
\se{dBV}\be\int_{t_0}^T\!\!\left(\|\partial{\bf V}\!/\partial t\|_s^2
+\|\partial{\bf B}/\partial t\|_s^2\right)^{\alpha_{s+2}/2}
\d t\le D_s^{(1)}\quad&\mbox{for~~}s\ge-1/2;\label{das}\\
\int_{t_0}^T\!\!\left(\|\partial{\bf V}\!/\partial t\|_s^2
+\|\partial{\bf B}/\partial t\|_s^2\right)^{\gamma_{s+5/2}/2}
\d t\le D_s^{(2)}\quad&\mbox{for~~}-5/2<s\le-1/2;\label{dgs}\\
\|\partial{\bf V}\!/\partial t\|_s^2+\|\partial{\bf B}/\partial t\|_s^2
\le D_s^{(3)}\quad&\mbox{for~~}s<-5/2;\label{ds}\\
\int_{t_0}^T\!\!\max_{\T}\Big(
\big|(-\nabla^2)^{s/2}\,\partial{\bf V}\!/\partial t\big|
+\big|(-\nabla^2)^{s/2}\,\partial{\bf B}/\partial t\big|
\Big)^{\!\!\alpha_{s+7/2}}\d t\le D_s^{(4)}\quad&\mbox{for~~}s>-2.\label{dtm}
\end{align}\end{subequations}
Here $D_s^{(i)}$ are sublinear functions of time $T\ge t_0$
that depend on the initial data and constants $s,\nu$ and $\eta$ only.

{\it Proof.} We have introduced a transformation \rf{svb} of coefficients
of the Fourier--Galerkin approximants $\bf V,\,B$ of solutions
to \rf{vb}. The modified coefficients satisfy equations \rf{GFm}, and
the time derivatives of $\bf V$ and $\bf B$ can be expanded as
\BE{\partial{\bf V}\over\partial t}=\sum_{0\ne|{\bf n}|\le N}
\ks^{\bf v}_{\bf n}(t)\e^{\i{\bf n\cdot x}-\delta\Phi|{\bf n}|},\qquad
{\partial{\bf B}\over\partial t}=\sum_{0\ne|{\bf n}|\le N}
\ks^{\bf b}_{\bf n}(t)\e^{{\i\bf n\cdot x}-\delta\Phi|{\bf n}|},\EE{eVB}
where it is denoted
\BE\ks^{\bf v}_{\bf n}(t)={\d\wv_{\bf n}\over\d t}
-\delta|{\bf n}|\wv_{\bf n}\,{\d\Phi\over\d t},\qquad
\ks^{\bf b}_{\bf n}(t)={\d\wb_{\bf n}\over\d t}
-\delta|{\bf n}|\wb_{\bf n}\,{\d\Phi\over\d t}.\EE{pd}

\subsection{Bounds in the Sobolev space norms}\label{bSn}

We need to bound Sobolev norms of the quantities
$$\ks^{\bf v}=\sum_{|{\bf n}|\le N}\ks^{\bf v}_{\bf n}\e^{\i\bf n\cdot x}
\mbox{~~and~~}
\ks^{\bf b}=\sum_{|{\bf n}|\le N}\ks^{\bf b}_{\bf n}\e^{\i\bf n\cdot x}.$$
Scalar multiplying \rf{mv} and \rf{mb} by $|{\bf n}|^{2s}\ks^{\bf v}_{\bf-n}$
and $|{\bf n}|^{2s}\ks^{\bf b}_{\bf-n}$, respectively, and summing up
the results over ${\bf n}\ne 0$ (note \rf{zm}) yields
\be&\|\ks^{\bf v}\|_s^2+\|\ks^{\bf b}\|_s^2
=-\sum_{\bf n}|{\bf n}|^{2s+2}(\nu\,\wv_{\bf n}\cdot\ks^{\bf v}_{\bf-n}
+\eta\,\wb_{\bf n}\cdot\ks^{\bf b}_{\bf-n})+\i\sum_{\bf n,k}|{\bf n}|^{2s}
\e^{\delta\Phi(|{\bf n}|-|{\bf k}|-|{\bf n-k}|)}\nonumber\\
&\hspace*{3ex}\times\!\left(\!
\big(\!\!-(\wv_{\bf n-k}\!\cdot\!{\bf k})\wv_{\bf k}
+(\wb_{\bf n-k}\!\cdot\!{\bf k})\wb_{\bf k}\big)\!\cdot\!\ks^{\bf v}_{\bf-n}
+\big((\wb_{\bf n-k}\!\cdot\!{\bf k})\wv_{\bf k}-(\wv_{\bf n-k}\!\cdot\!
{\bf k})\wb_{\bf k}\big)\!\cdot\!\ks^{\bf b}_{\bf-n}\right)\nonumber\\
\shortintertext{(by the orthogonality \rf{so})}
&\le{1\over4}\|\ks^{\bf v}\|_s^2+\nu^2\|\wv\|_{s+2}^2
+{1\over4}\|\ks^{\bf b}\|_s^2+\eta^2\|\wb\|_{s+2}^2\nonumber\\
&\hspace*{3ex}+\sum_{\bf n,k}\!|{\bf n}|^{2s+1-\beta}|{\bf k}|^\beta
\left(|\ks^{\bf v}_{\bf-n}|
\big(|\wv_{\bf k}||\wv_{\bf n-k}|+|\wb_{\bf k}||\wb_{\bf n-k}|\big)
+|\ks^{\bf b}_{\bf-n}|(|\wb_{\bf n-k}||\wv_{\bf k}|
+|\wv_{\bf n-k}||\wb_{\bf k}|)\right)\label{bb}\end{align}
for an arbitrary $\beta$ such that $0\le\beta\le 1$
(the triangle inequality is applied to bound the exponential, and the inequality
\rf{ob} is used with $\wv_{\bf n-k}$ and $\wb_{\bf n-k}$ replacing
$\bf v_{n-k}$ and $\bf b_{n-k}$). For different indices $s$ of the norms,
further derivations are similar, but differ in details.

{\it Proof of \rf{das}.} We assume $s=-1/2$ and $\beta=1/2$. By the embedding
theorem inequalities \rf{emb}, the last sum in \rf{bb} is majorised by
$$C'_{1/2}\!\left(\|\ks^{\bf v}\|_{-1/2}
\big(\|\wv\|_1^2+\|\wb\|_1^2\big)+2\|\ks^{\bf b}\|_{-1/2}
\|\wv\|_1\|\wb\|_1\!\right).$$
Consequently, \rf{bb} implies
\BE{1\over2}(\|\ks^{\bf v}\|_{-1/2}^2+\|\ks^{\bf b}\|_{-1/2}^2)
\le\nu^2\|\wv\|_{3/2}^2+\eta^2\|\wb\|_{3/2}^2+2(C'_{1/2})^2
(\|\wv\|_1^2+\|\wb\|_1^2)^2.\EE{mh}
By virtue of \rf{pd}, this inequality and \rf{exp},
\ba\|\partial{\bf V}\!/\partial t\|_s^2+\|\partial{\bf B}/\partial t\|_s^2
&=\sum_{{\bf n}\ne0}|{\bf n}|^{2s}\e^{-2\delta\Phi|{\bf n}|}
\big(|\ks^{\bf v}_{\bf n}|^2+|\ks^{\bf b}_{\bf n}|^2\big)\\
&\le(\|\ks^{\bf v}\|_{-1/2}^2+\|\ks^{\bf b}\|_{-1/2}^2)
\max_{{\bf n}\ne0}|{\bf n}|^{2s+1}\e^{-2\delta\Phi|{\bf n}|}\\
&\le2\left(\!{2s+1\over2\e\delta\Phi}\right)^{\!\!2s+1}\!\!\!
\big(\!\max(\nu^2,\eta^2)(\|\wv\|_{3/2}^2+\|\wb\|_{3/2}^2)
+2(C'_{1/2})^2(\|\wv\|_1^2+\|\wb\|_1^2)^2\big).\end{align*}
Therefore, by H\"older's inequality, \rf{ab}, \rf{iti} upon changing
$s\to s/2+5/4$, and \rf{Qp}
\ba&\int_{t_0}^T\!\!\big(\|\partial{\bf V}\!/\partial t\|_s^2
+\|\partial{\bf B}/\partial t\|_s^2\,\big)^{\!\alpha_{s+2}/2}\d t\\
&\hspace*{3ex}\le\!\int_{t_0}^T\!\!\Big(\!D_s^{(1,1)}\Big(
\big(\|\wv\|_{3/2}^2\!+\|\wb\|_{3/2}^2\big)^{\!{\alpha_{s+2}\over2}}\!+
\big(\|\wv\|_{3/2}^2\!+\|\wb\|_{3/2}^2\big)^{\!1/2}\Big)+D_s^{(1,2)}\,
\Phi^{-{2s+1\over2s+3}}\,(\|\wv\|_1^2\!+\|\wb\|_1^2)^{2\over2s+3}\Big)\d t\\
&\hspace*{3ex}\le D_s^{(1)}=D_s^{(1,1)}\Big((Q'(T))^{\alpha_{s+2}/2}
(T-t_0)^{1-\alpha_{s+2}/2}+Q'(T)\Big)+D_s^{(1,2)}Q''_{s/2+5/4}(T)
\end{align*}
for $s\ge-1/2$, where
$$D_s^{(1,1)}=\Big(\!2\max(\nu^2,\eta^2)\Big({2s+1\over2\e\delta}
\Big)^{\!2s+1}\Big)^{\!1/(2s+3)}\!,\qquad
D_s^{(1,2)}=(2C'_{1/2})^{2/(2s+3)}\Big({2s+1\over2\e\delta}
\Big)^{\!(2s+1)/(2s+3)}\!.$$

{\it Proof of \rf{dgs}.}
For $s$ in the subinterval $-1\!<\!s\!\le\!-1/2$, we assume $\beta\!=\!1\!+\!s$
and bound the last sum in \rf{bb} by
\se{mds}\be(2\pi)^{-3}&\int_{\T}\Big(f^{\ks^{\bf v}}_s\big(f^{\wv}_{1+s}f^{\wv}_0
+f^{\wb}_{1+s}f^{\wb}_0\big)+f^{\ks^{\bf b}}_s
(f^{\wv}_{1+s}f^{\wb}_0+f^{\wb}_{1+s}f^{\wv}_0\big)\Big)\d{\bf x}\nonumber\\
&\hspace*{-3em}\le\!\Big(\!\big|f^{\ks^{\bf v}}_s\big|_2
\big(\big|f^{\wv}_{1+s}\big|_{12\over5+2s}\big|f^{\wv}_0\big|_{12\over1-2s}
\!+\big|f^{\wb}_{1+s}\big|_{12\over5+2s}\big|f^{\wb}_0\big|_{12\over1-2s}\big)
+\big|f^{\ks^{\bf b}}_s\big|_2
\big(\big|f^{\wv}_{1+s}\big|_{12\over5+2s}\big|f^{\wb}_0\big|_{12\over1-2s}
\!+\big|f^{\wb}_{1+s}\big|_{12\over5+2s}\big|f^{\wv}_0\big|_{12\over1-2s}\big)\!\Big)\nonumber\\
\shortintertext{(see \rf{scf})}
&\le\widetilde C'_s\Big(\|\ks^{\bf v}\|_s
\big(\|\wv\|_{5+2s\over4}^2+\|\wb\|_{5+2s\over4}^2\big)
+2\|\ks^{\bf b}\|_s\|\wv\|_{5+2s\over4}\|\wb\|_{5+2s\over4}\Big),\label{ms1}\\
\shortintertext{where the constant
$\widetilde C'_s=C_{(1-2s)/4}C_{(2s+5)/4}$ has been introduced.
For $-5/2<s\le-1$, we assume $\beta=0$ and majorise
the last sum in \rf{bb} as follows:}
(2\pi)^{-3}&\int_{\T}\Big(f^{\ks^{\bf v}}_{2s+1}\big((f^{\wv}_0)^2
+(f^{\wb}_0)^2\big)+2f^{\ks^{\bf b}}_{2s+1}f^{\wv}_0f^{\wb}_0\Big)\d{\bf x}\nonumber\\
&\le\Big(\big|f^{\ks^{\bf v}}_{2s+1}\big|_{6\over5+2s}\big(
\big|f^{\wv}_0\big|^2_{12\over1-2s}+\big|f^{\wb}_0\big|^2_{12\over1-2s}\big)
+2\big|f^{\ks^{\bf b}}_{2s+1}\big|_{6\over5+2s}
\big|f^{\wv}_0\big|_{12\over1-2s}\big|f^{\wb}_0\big|_{12\over1-2s}\Big)\nonumber\\
&\le\widetilde C'_s\Big(\|\ks^{\bf v}\|_s
\big(\|\wv\|_{5+2s\over4}^2+\|\wb\|_{5+2s\over4}^2\big)
+2\|\ks^{\bf b}\|_s\|\wv\|_{5+2s\over4}\|\wb\|_{5+2s\over4}\Big),
\label{ms2}\end{align}\end{subequations}
where $\widetilde C'_s=C_{-1-s}C_{(2s+5)/4}^2$.
Applying now \rf{ms1} or \rf{ms2} (depending on to which of the two
subintervals $s$ belongs), we infer from \rf{bb}
\be\|\ks^{\bf v}\|_s^2+\|\ks^{\bf b}\|_s^2
&\le2\nu^2\|\wv\|_{s+2}^2+2\eta^2\|\wb\|_{s+2}^2
+4(\widetilde C'_s)^2\left(\|\wv\|_{2s+5\over4}^2+\|\wb\|_{2s+5\over4}^2
\right)^{\!\!2}\nonumber\\
&\le2\nu^2\|\wv\|_{s+2}^2+2\eta^2\|\wb\|_{s+2}^2
+4(\widetilde C'_s)^2\left(\|\wv\|_1^{2s+5\over2}\|\wv\|_0^{-2s-1\over2}
\!\!+\|\wb\|_1^{2s+5\over2}\|\wb\|_0^{-2s-1\over2}\right)^{\!\!2}
\label{dV}\end{align}
(by H\"older's inequality) for all the considered $s$ in the interval
$-5/2<s\le-1/2$.

For $-5/2<s\le-1/2$, the last term (arising from the nonlinear terms
in \rf{vb}) in the r.h.s.~of \rf{dV}
becomes time-integrable upon raising to the power $\gamma_{s+5/2}/2=2/(2s+5)$.
For $-1\le s\le-1/2$, the other two terms in the r.h.s.~of \rf{dV},
arising from the linear diffusivity terms in \rf{vb}, are time-integrable
when raised to the higher power $\alpha_{s+2}/2=1/(2s+3)\ge2/(2s+5)$.
For $-2<s\le-1$,
$$\nu^2\|\wv\|_{s+2}^2+\eta^2\|\wb\|_{s+2}^2\le
\nu^2\|\wv\|_1^{2(s+2)}\|\wv\|_0^{-2(s+1)}
+\eta^2\|\wb\|_1^{2(s+2)}\|\wb\|_0^{-2(s+1)},$$
and therefore the terms describing diffusivity become integrable in time
if raised to the power $1/(s+2)\ge2/(2s+5)$; for $s\le-2$ they are finite
at any time due to the energy inequality.
Thus, for $s\le-1/2$ they do not affect the maximum power to which
the l.h.s.~of \rf{dV} can be raised without losing the time integrability. We
conclude, in view of relations \rf{eVB}, that for $-5/2<s\le-1/2$ the integrals
$$\int_{t_0}^T\!\!\Big(\|{\partial{\bf V}\!/\partial t\|_s^2
+\|\partial{\bf B}/\partial t}\|_s^2\Big)^{\!\!\gamma_{s+5/2}/2}\d t
\le\int_{t_0}^T\!\big(\|\ks^{\bf v}\|_s^2+\|\ks^{\bf b}\|_s^2\big)^{\!\gamma_{s+5/2}/2}\d t$$
remain bounded for all $T>t_0$; it is easy to deduce explicit expressions
for $D_s^{(2)}$ in \rf{dgs} from \rf{dV} and \rf{bou}. Clearly, for $s=-5/2$
the integrals remain finite when raised to any positive power.

{\it Proof of \rf{ds}.}
We assume $s<-5/2$ and $\beta=0$. For $p>3/2$, by the embedding theorem,
\BE\max_{{\bf x}\in\T}|f({\bf x})|<c_p\|f\|_p\qquad\mbox{for any~~}
f\in H_p(\T),\EE{eT}
where $c_p$ is a constant that depends on $p$ only. Applying this inequality
to the last sum in \rf{bb}, we find its upper bound
\ba\max_{\T}&|f^{\ks^{\bf v}}_{2s+1}|\,(2\pi)^{-3}\int_{\T}
\big((f^{\wv}_0)^2+(f^{\wb}_0)^2\big)\d{\bf x}
+2\max_{\T}|f^{\ks^{\bf b}}_{2s+1}|\,(2\pi)^{-3}\int_{\T}
|f^{\wb}_0f^{\wv}_0|\d{\bf x}\\
&\le c_{-s-1}\big(\|\ks^{\bf v}\|_s
+\|\ks^{\bf b}\|_s\big)(\|\wv\|_0^2+\|\wb\|_0^2),\end{align*}
whereby \rf{bb} implies (since the terms in the r.h.s.~of \rf{bb}, related to
diffusivity, have a negative norm index $s+1$)
$$\|\ks^{\bf v}\|_s^2+\|\ks^{\bf b}\|_s^2
\le D_s^{(3)}=2\Big(\!\!\max(\nu^2,\eta^2)\big(\|\wv\|_0^2+\|\wb\|_0^2\big)
+c_{-s-1}^2\big(\|\wv\|_0^2+\|\wb\|_0^2\big)^2\Big).$$

\subsection{Bounds in the Wiener algebra norms}\label{dWn}

{\it Proof of \rf{dtm}.}
Once again, we bound the maxima by the sums of absolute values of the Fourier
coefficients of the respective functions. By the Lemma and \rf{mh}, for $s>-2$,
\ba|(-\nabla^2)^{s/2}\,\partial{\bf V}\!/\partial t|
+|(-\nabla^2)^{s/2}\,\partial{\bf B}/\partial t|
&=\Big|\sum_{{\bf n}\ne0}|{\bf n}|^s
\ks^{\bf v}_{\bf n}\e^{-\delta\Phi|{\bf n}|+\i\bf n\cdot x}\Big|
+\Big|\sum_{{\bf n}\ne0}|{\bf n}|^s
\ks^{\bf b}_{\bf n}\e^{-\delta\Phi|{\bf n}|+\i\bf n\cdot x}\Big|\\
&\le\Big(\!\Big(\sum_{{\bf n}\ne0}{|\ks^{\bf v}_{\bf n}|^2\over|{\bf n}|}
\Big)^{\!\!1/2}\!+\Big(\sum_{{\bf n}\ne0}{|\ks^{\bf b}_{\bf n}|^2\over|{\bf n}|}
\Big)^{\!\!1/2}\Big)\!\Big(\sum_{{\bf n}\ne0}|{\bf n}|^{2s+1}
\e^{-2\delta\Phi|{\bf n}|}\Big)^{\!\!1/2}\\
&\le\sqrt2\big(\|\ks^{\bf v}_{\bf n}\|^2_{-1/2}
+\|\ks^{\bf b}_{\bf n}\|^2_{-1/2}\big)^{\!1/2}C_{2s+1,2\delta}\,\Phi^{-(s+2)}.\end{align*}
Thus, for $s>-2$, by \rf{mh}, \rf{ab}, H\"older's inequality,
\rf{iti} upon changing $s\to s/2+2$, and \rf{Qp},
\ba\int_{t_0}^T&\!\!\max_{\T}\Big(
|(-\nabla^2)^{s/2}\,\partial{\bf V}\!/\partial t|
+|(-\nabla^2)^{s/2}\,\partial{\bf B}/\partial t|\Big)^{\!\!\alpha_{s+7/2}}\d t\\
&\le\int_{t_0}^T\!\!\Big(\!D_s^{(4,1)}
\big((\|\wv\|_{3/2}^2\!+\|\wb\|_{3/2}^2)^{1\over s+3}
+(\|\wv\|_{3/2}^2\!+\|\wb\|_{3/2}^2)^{1/2}\big)
+\!D_s^{(4,2)}(\|\wv\|_1^2\!+\|\wb\|_1^2)^{1\over s+3}\Phi^{-{s+2\over s+3}}
\Big)\d t\\
&\le D_s^{(4)}=D_s^{(4,1)}\Big((Q'(T))^{1\over s+3}(T-t_0)^{s+2\over s+3}
+Q'(T)\Big)+D_s^{(4,2)}Q''_{s/2+2}(T),\end{align*}
where $D_s^{(4,1)}=(2C_{2s+1,2\delta}\max(\nu,\eta))^{\alpha_{s+7/2}}$ and
$D_s^{(4,2)}=(\sqrt8\,C_{2s+1,2\delta}C'_{1/2})^{\alpha_{s+7/2}}.~~\RHD$

\subsection{Bounds for time derivatives stemming
from the inequalities \rf{qq}}\label{alt}

A priori bounds for higher-index norms of $\d{\bf V}\!/\d t$ and $\d{\bf B}/\d t$ can also
be constructed following the standard techniques by using \rf{vb} directly.
We show here that for $s\ge-1/2$ this yields inequalities similar to \rf{das}
and involving the same exponents $\alpha_{s+2}$ (for which the $H_s(\T)$ norms
of the second derivatives are guaranteed to be time-integrable).

We scalar multiply \rf{wV} and \rf{wB} by $|{\bf n}|^{2s}\d\hV_{\bf-n}/\d t$
and $|{\bf n}|^{2s}\d\hB_{\bf-n}/\d t$, respectively, and sum up the results
over ${\bf n}\ne 0$ (by virtue of \rf{zm}). Let us denote
$$\widetilde C_s=\left\{\begin{array}{ll}
1&\mbox{~~for~~}-1\le s\le0;\\2^{s+1}&\mbox{~~for~~}s>0.\end{array}\right.$$
Taking into account the orthogonality \rf{so}, the inequalities
$$|{\bf n}|^{2s}|\hB_{\bf n-k}\!\cdot\!{\bf k}|\le\left\{\begin{array}{ll}
|{\bf n}|^s|{\bf k}|^{1+s}|\hB_{\bf n-k}|&\mbox{~~for~~}-1\le s\le0;\\
2^s|{\bf n}|^s(|{\bf k}|^{s+1}+|{\bf n-k}|^{s+1})
|\hB_{\bf n-k}|&\mbox{~~for~~}s>0\end{array}\right.$$
(stemming from \rf{ob}, where $\beta=s$ and $\hV_{\bf n-k}$ and $\hB_{\bf n-k}$
replace $\bf v_{n-k}$ and $\bf b_{n-k}$, and from \rf{pin}) and similar ones
for $\hV_{\bf n-k}$,
and the embedding theorem inequalities \rf{emb}, we obtain for $s\ge-1$
\be&\left\|{\partial{\bf V}\over\partial t}\right\|_s^2
+\left\|{\partial{\bf B}\over\partial t}\right\|_s^2
=-\sum_{\bf n}|{\bf n}|^{2s+2}
\left(\nu\hV_{\bf n}\cdot{\d\hV_{\bf-n}\over\d t}
+\eta\hB_{\bf n}\cdot{\d\hB_{\bf-n}\over\d t}\right)\nonumber\\
&\hspace*{3ex}+\i\sum_{\bf n,k}|{\bf n}|^{2s}\!\left(\!\!
\Big(\!\!-(\hV_{\bf n-k}\!\cdot\!{\bf k})\hV_{\bf k}
+(\hB_{\bf n-k}\!\cdot\!{\bf k})\hB_{\bf k}\Big)\!\cdot\!
{\d\hV_{\bf-n}\over\d t}
+\Big(\!(\hB_{\bf n-k}\!\cdot\!{\bf k})\hV_{\bf k}
-(\hV_{\bf n-k}\!\cdot\!{\bf k})\hB_{\bf k}\Big)\!\cdot\!
{\d\hB_{\bf-n}\over\d t}\right)\nonumber\\
&\le{1\over4}\left\|{\partial{\bf V}\over\partial t}\right\|_s^2
+\nu^2\|{\bf V}\|_{s+2}^2+{1\over4}\left\|{\partial{\bf B}\over
\partial t}\right\|_s^2+\eta^2\|{\bf B}\|_{s+2}^2\nonumber\\
&\hspace*{3ex}+\widetilde C_s\!\sum_{\bf n,k}|{\bf n}|^s|{\bf k}|^{s+1}\!
\left(\vphantom{\partial\over\partial}\right|\!
{\partial{\hV_{-n}}\over\partial t}\!
\left|\vphantom{\partial\over\partial}\right.
\!\Big(|\hV_{\bf n-k}||\hV_{\bf k}|+|\hB_{\bf n-k}||\hB_{\bf k}|\Big)
\!+\left|\vphantom{\partial\over\partial}\right.\!\!
{\partial{\hB_{-n}}\over\partial t}\!
\left|\vphantom{\partial\over\partial}\right.
\!\Big(|\hB_{\bf k}||\hV_{\bf n-k}|+|\hV_{\bf k}||\hB_{\bf n-k}|\Big)\!\!\!
\left)\vphantom{\partial\over\partial}\right.\nonumber\\
\shortintertext{(for $s>0$, we have used the invariance of the last sum
under the change of the index $\bf k\to n-k$)}
&\le{1\over4}\left\|{\partial{\bf V}\over\partial t}\right\|_s^2
+\nu^2\|{\bf V}\|_{s+2}^2+{1\over4}\left\|{\partial{\bf B}\over\partial t}\right\|_s^2
+\eta^2\|{\bf B}\|_{s+2}^2\nonumber\\
&\hspace*{3ex}+C'_{1/2}\widetilde C_s\!\left(\left\|{\partial{\bf V}\over\partial t}\right\|_s
\!\!\!\big(\|{\bf V}\|_{s+3/2}\|{\bf V}\|_1+\|{\bf B}\|_{s+3/2}\|{\bf B}\|_1\big)
+\!\left\|{\partial{\bf B}\over\partial t}\right\|_s
\!\!\!\big(\|{\bf B}\|_{s+3/2}\|{\bf V}\|_1+\|{\bf V}\|_{s+3/2}\|{\bf B}\|_1\big)
\!\!\right)\nonumber\\
&\le{1\over2}\left\|{\partial{\bf V}\over\partial t}\right\|_s^2
\!\!+\nu^2\|{\bf V}\|_{s+2}^2\!+{1\over2}\left\|{\partial{\bf B}\over\partial t}\right\|_s^2
\!\!+\eta^2\|{\bf B}\|_{s+2}^2\!+{C'''_s\over2}
(\|{\bf V}\|_{s+3/2}^2\!+\|{\bf B}\|_{s+3/2}^2)
(\|{\bf V}\|_1^2\!+\|{\bf B}\|_1^2),\label{mi}\end{align}
where $C'''_s=(2C'_{1/2}\widetilde C_s)^2$. Now Young's inequality,
the identity $1-\alpha_{s+p}/\alpha_{s+q}=(p-q)\alpha_{s+p}$ and
the inequality \rf{ab} yield for $s\ge-1/2$
\ba&\hspace*{-3ex}\left(\left\|{\partial{\bf V}\over\partial t}\right\|_s^2
+\left\|{\partial{\bf B}\over\partial t}\right\|_s^2\right)^{\!\alpha_{s+2}/2}
\le(\sqrt2\max(\nu,\eta))^{\alpha_{s+2}}
(\|{\bf V}\|_{s+2}^2+\|{\bf B}\|_{s+2}^2)^{\alpha_{s+2}/2}\\
&+(C'''_s\alpha_{s+2}/\alpha_{s+3/2})^{\alpha_{s+2}/2}
(\|{\bf V}\|_{s+3/2}^2+\|{\bf B}\|_{s+3/2}^2)^{\alpha_{s+3/2}/2}
\!+(C'''_s\alpha_{s+2}/2)^{\alpha_{s+2}/2}(\|{\bf V}\|_1^2+\|{\bf B}\|_1^2),
\end{align*}
implying a bound of the form \rf{das}.

\section{From a priori bounds to bounds for weak solutions}\label{bw}

The goal of this section is to prove that a priori bounds \rf{dBV}
for the time derivatives of the Fourier--Galerkin approximants
of solutions to the problem \rf{vb} are also satisfied by the derivatives
of the weak solutions. In this section, the dependence of the approximants
on the resolution parameter $N$ is shown explicitly.

\subsection{Justification of the bounds \rf{zz} for the weak solutions}

It is instructive to recall how weak solutions to \rf{vb} are constructed.

{\it Theorem 5}. For any $T>0$, there exists a subsequence $N_j\to\infty$ such
that, for all $\bf n$, $\hV^{(N_j)}_{\bf n}$ and $\hB^{(N_j)}_{\bf n}$ converge
pointwise uniformly on $[0,T]$ to continuous functions $\hV_{\bf n}$ and
$\hB_{\bf n}$, respectively. The fields \rf{Fs}
are weak solutions to equations \rf{vb}. They belong to $L_2(\T)$ at any time
and are weakly continuous in time in $L_2(\T)$;
$\partial{\bf V}\!/\partial x_m,\partial{\bf B}/\partial x_m\in
L_2(\T\times[0,T])$ for all $m$. The energy inequality
\BE{1\over2}(\|{\bf V}\|_0^2+\|{\bf B}\|_0^2)+\int_0^T\big(
\nu\|{\bf V}\|_1^2+\eta\|{\bf B}\|_1^2\big)\,\d t\le
E={1\over2}(\|{\bf V}^{\rm(init)}\|_0^2+\|{\bf B}^{\rm(init)}\|_0^2)\EE{en}
is satisfied, as well as the bounds \rf{zz}.

{\it Proof.} Coefficients of the truncated series \rf{Fou} satisfy equations
\rf{wF} which we consider on a time interval $0\le t\le T$.
Scalar multiplying \rf{wV} and \rf{wB} by $2\hV_{\bf-n}^{(N)}$ and
$2\hB_{\bf-n}^{(N)}$, respectively, and summing up yields
\BE\|{\bf V}^{(N)}\|_0^2+\|{\bf B}^{(N)}\|_0^2+2\int_0^T\!\!\!\!\big(
\nu\|{\bf V}^{(N)}\|_1^2+\eta\|{\bf B}^{(N)}\|_1^2\big)\,\d t\!
=\!\|{\bf V}^{{\rm(init},N)}\|_0^2+\|{\bf B}^{{\rm(init},N)}\|_0^2
\!\le\!2E.\EE{en0}
Here ${\bf V}^{{\rm(init},N)}$ and ${\bf B}^{{\rm(init},N)}$ denote
the initial fields upon projecting them onto the subspace, in which
the Fourier--Galerkin approximant of the solution is sought.
In view of this inequality, integrating \rf{wV} over a time interval
$t_1\le t\le t_2$ such that $0\le t_1<t_2\le T$ yields, for $N\ge|\bf n|$,
\se{UB}\be|\hV_{\bf n}^{(N)}(t_1)-\hV_{\bf n}^{(N)}(t_2)|\le&
\int_{t_1}^{t_2}\!\!\Big(\nu|{\bf n}|^2|\hV_{\bf n}^{(N)}|
+\!\sum_{\bf k}\big|(\hV_{\bf k}^{(N)}\!\cdot\!{\bf n})\hV_{\bf n-k}^{(N)}
-(\hB_{\bf k}^{(N)}\!\cdot\!{\bf n})\hB_{\bf n-k}^{(N)}\big|\Big)\d t\nonumber\\
\le&\,|t_1-t_2|\Big(\nu|{\bf n}|^2\max_{t_1\le t\le t_2}|\hV_{\bf n}^{(N)}|
+|{\bf n}|\max_{t_1\le t\le t_2}\big(\|\hV^{(N)}\|_0^2
+\|\hB^{(N)}\|_0^2\big)\!\Big)\nonumber\\
\le&\,|t_1-t_2|\big(\nu|{\bf n}|^2\sqrt{2E}+2E|{\bf n}|\big);\label{Vub}
\shortintertext{similarly, from \rf{wB},}
|\hB_{\bf n}^{(N)}(t_1)-\hB_{\bf n}^{(N)}(t_2)|\le&
\,|t_1-t_2|\big(\eta|{\bf n}|^2\sqrt{2E}
+E|{\bf n}|\big).\label{Bub}\end{align}\end{subequations}
Inequality \rf{en0} implies that for each wave vector $\bf n$ coefficients
${\bf V}_{\bf n}^{(N)}$, as well as ${\bf B}_{\bf n}^{(N)}$, are uniformly
(over the resolution $N$) bounded on the given time interval.
For each wave vector $\bf n$, by \rf{UB}, the functional
sets $\big\{\hV_{\bf n}^{(N)}\big|N\ge|\bf n|\big\}$
and $\big\{\hB_{\bf n}^{(N)}\big|N\ge|\bf n|\big\}$
are equicontinuous uniformly over $N$. Apply\-ing the Arzel\`a--Ascoli theorem
and the diagonal process, we can extract a subsequence $N_j\to\infty$
such that, for any $\bf n$, the approximants $\hV_{\bf n}^{(N_j)}$ uniformly
on $[0,T]$ converge to a continuous in time limit function $\hV_{\bf n}$,
and similarly $\hB_{\bf n}^{(N_j)}\to\hB_{\bf n}$. Now \rf{Fs} are weak
solutions to \rf{vb}.

They must satisfy integral identities obtained by scalar multiplying \rf{ve}
and \rf{be} in $L^2(\T)$ by arbitrary smooth space-periodic solenoidal
test vector fields $\bf f_V$ and $\bf f_B$, respectively, and~integrating by parts
over the cylinder $\T\times[0,T]$ so that the test fields only would be
differentiated in the integrand. These identities can be proven by the standard
arguments, taking the limit $N_j\to\infty$ in \rf{wF} and recalling that
convergence ${\bf V}_{\bf n}^{(N_j)}\to\bf V$ and
${\bf B}_{\bf n}^{(N_j)}\to\bf B_n$ is uniform in time, and the embedding
$H_1(\T)\to L_2(\T)$ is compact. We do not give a detailed proof here.

The equality \rf{en0} does not necessarily hold for weak solutions
to \rf{vb}, but it implies the inequality \rf{en}.
To see this, we consider partial sums truncated at a certain level $M\le N$:
$$\sum_{|{\bf n}|\le M}\big(|{\bf V}_{\bf n}^{(N)}|^2
+|{\bf B}_{\bf n}^{(N)}|^2\big)+2\int_0^T\sum_{|{\bf n}|\le M}|{\bf n}|^2
\big(\nu|{\bf V}_{\bf n}^{(N)}|^2+\eta|{\bf B}_{\bf n}^{(N)}|^2\big)\,\d t\le2E.$$
In the limit $N_j\to\infty$ for the chosen subsequence, this inequality
takes the form
$$\sum_{|{\bf n}|\le M}\big(|{\bf V}_{\bf n}|^2+|{\bf B}_{\bf n}|^2\big)
+2\int_0^T\sum_{|{\bf n}|\le M}|{\bf n}|^2
\big(\nu|{\bf V}_{\bf n}|^2+\eta|{\bf B}_{\bf n}|^2\big)\,\d t\le2E.$$
Since the latter inequality holds true for all truncation parameter values $M$,
we obtain \rf{en}, whereby ${\bf V(x},t),{\bf B(x},t)\in L_2(\T)$ at any time,
and $\partial{\bf V}\!/\partial x_m,\partial{\bf B}/\partial x_m$ belong to
$L_2(\T\times[0,T])$.

To establish weak continuity of ${\bf V(x},t)$ in time at time $t$, it suffices
to show that, given a field ${\bf a}\in L_2(\T)$, we can find
$\tau(t,{\bf a})>0$ such that
$$\zeta(\tau)=\Big|\int_{\T}\big({\bf V}(t+\tau)
-{\bf V}(t)\big)\cdot{\bf a}\,\d{\bf x}\,\Big|$$
is below any given threshold. We split
${\bf a}=\sum_{|{\bf n}|\le M}\widehat{\bf a}_{\bf n}\e^{\i\bf n\cdot x}
+{\bf a}'$. By \rf{Vub},
\ba\zeta(\tau)&\le(2\pi)^3\Big|\sum_{|{\bf n}|\le M}\!\big(\hV_{\bf n}(t+\tau)
-\hV_{\bf n}(t)\big)\cdot\widehat{\bf a}_{\bf-n}\Big|+\Big|\int_{\T}\!
\big({\bf V}(t+\tau)-{\bf V}(t)\big)\cdot{\bf a}'\d{\bf x}\Big|\\
&\le\tau(2\pi)^3\sum_{|{\bf n}|\le M}\!
\big(\nu|{\bf n}|^2\sqrt{2E}+E|{\bf n}|\big)|\widehat{\bf a}_{\bf-n}|
+2(2\pi)^3\|{\bf a}'\|_0\sqrt{2E}.\end{align*}
On increasing $M$, $\|{\bf a}'\|_0$ becomes sufficiently small,
and for this $M$ the first term is made sufficiently small by choosing
an appropriate $\tau$. Weak continuity of $\bf B$ is established the same way.

Like the energy inequality, \rf{zz} can be proven for the weak solution
by passing to the limit $N_j\to\infty$ in the a priori inequalities \rf{zz}
for the approximants, where the norms of the approximants are replaced
by the respective sums over wave vectors for $|{\bf n}|\le M$. In the case
of \rf{was}, we apply this procedure to the stronger inequality, where
sums of absolute values of the Fourier coefficients of the respective terms
replace the maxima in the l.h.s.~~$\RHD$

Furthermore, suppose $\bf V$ and $\bf B$ belong to $\dot H_s(\T)$ at time
$t=t_0$ for some $s$ such that $1/2<s\le1$. Then \rf{tsu}, \rf{hsb} and
\rf{exp} imply that for $t<t_0+t_*(t_0)$ and any $p\ge s$
$$\|{\bf V}^{(N)}\|_p^2+\|{\bf B}^{(N)}\|_p^2
\le(\|{\bf v}\|_s^2+\|{\bf b}\|_s^2)\max_{\bf n}|{\bf n}|^{2(p-s)}
\e^{-2\sigma|{\bf n}|(t-t_0)}\le q_s(t-t_0)((p-s)/(\e\sigma(t-t_0)))^{2(p-s)},$$
where
\BE t_*(t_0)=(\|{\bf V}(t_0)\|_s^2+\|{\bf B}(t_0)\|_s^2)^{-2/(2s-1)}/C''_s>0\EE{tx}
(cf.~\rf{tb}). It is legitimate to pass to the limit $N_j\to\infty$ to obtain
\BE\|{\bf V}\|_p^2+\|{\bf B}\|_p^2\le q_s(t-t_0)\,\big((p-s)/(\e\sigma(t-t_0))
\big)^{2(p-s)}=n_p(t-t_0;t_0)\quad\mbox{for~~}t_0<t<t_0+t_*(t_0)\EE{qu}
(the second argument in $n_p$ reflects that $q_s$ involves norms
of the solution at $t=t_0$).

\subsection{A bound for $\|\partial^2{\bf V}^{(N)}\!/\partial t^2\|_s$ and
$\|\partial^2{\bf B}^{(N)}\!/\partial t^2\|_s$ for $s<-7/2$}

We may try to apply a similar reasoning to establish bounds \rf{dBV} for time
derivatives of a weak solution to \rf{vb}. By \rf{ds}, for each wave vector
$\bf n$, the derivatives $\d{\bf V}_{\bf n}^{(N)}\!/\d t$ and
$\d{\bf B}_{\bf n}^{(N)}\!/\d t$ are uniformly (over the resolution parameter
$N$) bounded on $[0,T]$. We could apply
the Arzel\`a--Ascoli theorem, if we showed that for each wave vector $\bf n$
the functional sets $\{\d\hV_{\bf n}^{(N)}\!/\d t\}$ and
$\{\d\hB_{\bf n}^{(N)}\!/\d t\}$ are uniformly (over~$N\ge|\bf n|$)
equicontinuous. We need, therefore, bounds for second time derivatives
of the approximants. Differentiating \rf{wF} in time yields
\se{dGF}\be{\d^2\hV_{\bf n}^{(N)}\over\d t^2}=&
-\nu|{\bf n}|^2{\d\hV_{\bf n}^{(N)}\over\d t}-\i\sum_{{\bf k}\ne0}\P\!\left(\!\!
\left(\!{\d\hV_{\bf n-k}^{(N)}\over\d t}\!\cdot\!{\bf n}\right)\!\hV_{\bf k}^{(N)}
+\left(\!\hV_{\bf n-k}^{(N)}\!\cdot\!{\bf n}\!\right){\d\hV_{\bf k}^{(N)}\over\d t}
\right.\nonumber\\
&\left.-\left(\!{\d\hB_{\bf n-k}^{(N)}\over\d t}\!\cdot\!{\bf n}\!\right)\!\hB_{\bf k}^{(N)}
-\left(\hB_{\bf n-k}^{(N)}\!\cdot\!{\bf n}\right)\!{\d\hB_{\bf n}^{(N)}\over\d t}
\!\right)\!,\label{dFv}\\
{\d^2\hB_{\bf n}^{(N)}\over\d t^2}=&
-\eta|{\bf n}|^2{\d\hB_{\bf n}^{(N)}\over\d t}
+\i{\bf n}\times\!\sum_{{\bf k}\ne0}\left(\!
{\d\hV_{\bf n-k}^{(N)}\over\d t}\!\times\hB_{\bf k}^{(N)}
+\hV_{\bf k}^{(N)}\!\times{\d\hB_{\bf n-k}^{(N)}\over\d t}\right)\!.
\label{dFb}\end{align}\end{subequations}

The r.h.s.~of \rf{dGF} involve first derivatives of the Fourier coefficients.
Their bounds in the space $\dot H_{-1}(\T)$ fit best our goals.
We obtain from \rf{mi}
\BE\|\partial{\bf V}^{(N)}\!/\partial t\|_{-1}^2
+\|\partial{\bf B}^{(N)}\!/\partial t\|_{-1}^2
\le\widetilde C''\big(\|{\bf V}^{(N)}\|_1^2+\|{\bf B}^{(N)}\|_1^2
+(\|{\bf V}^{(N)}\|_1^2+\|{\bf B}^{(N)}\|_1^2)^{3/2}\big),\EE{clo}
where $\widetilde C''$ depends only on the parameters
of the problem (the diffusivities $\nu$ and~$\eta$) and the initial data
${\bf V}^{\rm(init)}$ and ${\bf B}^{\rm(init)}$.

Bounds for the norms $\|\partial^2{\bf V}^{(N)}\!/\partial t^2\|_s$ and
$\|\partial^2{\bf B}^{(N)}\!/\partial t^2\|_s$ of any index are suitable
to establish equicontinuity for a fixed wave vector $\bf n$.
We choose for simplicity $s<-7/2$, because this gives an opportunity to employ
the embedding theorem inequality \rf{eT}. We scalar multiply \rf{dFv}
and \rf{dFb} by $|{\bf n}|^{2s}\d^2\hV_{\bf-n}^{(N)}/\d t^2$ and
$|{\bf n}|^{2s}\d^2\hB_{\bf-n}^{(N)}/\d t^2$, respectively, use the inequality
$2|{\bf n}||{\bf k}|\ge|{\bf n-k}|$ valid for $|{\bf n}|\ge1$ and
$|{\bf k}|\ge1$, sum up the results over ${\bf n}\ne 0$, and obtain
\ba&\left\|{\partial^2{\bf V}^{(N)}\over\partial t^2}\right\|_s^2
\!\!+\left\|{\partial^2{\bf B}^{(N)}\over\partial t^2}\right\|_s^2
\!\!+{1\over2}\,{\d\over\d t}\left(
\nu\left\|{\partial{\bf V}^{(N)}\over\partial t}\right\|_{s+1}^2
\!\!+\eta\left\|{\partial{\bf B}^{(N)}\over\partial t}\right\|_{s+1}^2
\right)\nonumber\\
&\hspace*{3em}\le4\left(\sum_{{\bf n}\ne0}|{\bf n}|^{2s+2}\,
\left|{\partial^2\hV_{\bf n}^{(N)}\over\partial t^2}\right|\right)\,
\max_{{\bf n}\ne0}\sum_{\bf k}{|{\bf k}|\over|{\bf n-k}|}\left(
\left|{\d\hV_{\bf n-k}^{(N)}\over\d t}\right|\!\big|\hV_{\bf k}^{(N)}\big|
+\left|{\d\hB_{\bf n-k}^{(N)}\over\d t}\right|\!\big|\hB_{\bf k}^{(N)}\big|
\right)\nonumber\\
&\hspace*{4em}+2\left(\sum_{{\bf n}\ne0}|{\bf n}|^{2s+2}\,
\left|{\partial^2\hB_{\bf n}^{(N)}\over\partial t^2}\right|\right)\,
\max_{{\bf n}\ne0}\sum_{\bf k}{|{\bf k}|\over|{\bf n-k}|}\left(
\left|{\d\hV_{\bf n-k}^{(N)}\over\d t}\right|\!\big|\hB_{\bf k}^{(N)}\big|
+\big|\hV_{\bf k}^{(N)}\big|\!\left|{\d\hB_{\bf n-k}^{(N)}\over\d t}\right|
\right)\\
&\hspace*{3em}\le
{1\over2}\left\|{\partial^2{\bf V}^{(N)}\over\partial t^2}\right\|_s^2
+8c^2_{-s-2}\left(\left\|{\partial{\bf V}^{(N)}\over\partial t}\right\|_{-1}
\!\!\|{\bf V}^{(N)}\|^{\phantom{2}}_1
+\left\|{\partial{\bf B}^{(N)}\over\partial t}\right\|_{-1}\!\!
\|{\bf B}^{(N)}\|^{\phantom{2}}_1\right)^{\!\!\!2}\nonumber\\
&\hspace*{4em}
+{1\over2}\left\|{\partial^2{\bf B}^{(N)}\over\partial t^2}\right\|_s^2
+2c_{-s-2}^2\left(\left\|{\partial{\bf V}^{(N)}\over\partial t}\right\|_{-1}\!\!
\|{\bf B}^{(N)}\|^{\phantom{2}}_1
+\left\|{\partial{\bf B}^{(N)}\over\partial t}\right\|_{-1}\!\!
\|{\bf V}^{(N)}\|^{\phantom{2}}_1\right)^{\!\!\!2}\end{align*}
(by \rf{eT}), and thus, by \rf{clo},
\be&\|\partial^2{\bf V}^{(N)}\!/\partial t^2\|_s^2
+\|\partial^2{\bf B}^{(N)}\!/\partial t^2\|_s^2
+{\d\over\d t}\Big(\nu\|\partial{\bf V}^{(N)}\!/\partial t\|_{s+1}^2
\!+\eta\|\partial{\bf B}^{(N)}\!/\partial t\|_{s+1}^2\Big)
\nonumber\\
&\hspace*{4em}\le20c^2_{-s-2}\left(
\left\|\partial{\bf V}^{(N)}\!/\partial t\right\|_{-1}^2
+\left\|\partial{\bf B}^{(N)}\!/\partial t\right\|_{-1}^2\right)
\big(\|{\bf V}^{(N)}\|_1^2+\|{\bf B}^{(N)}\|_1^2\big)\nonumber\\
&\hspace*{4em}\le20c^2_{-s-2}\,
\widetilde C''\Big(\!\big(\|{\bf V}^{(N)}\|_1^2+\|{\bf B}^{(N)}\|_1^2\big)^2
+\big(\|{\bf V}^{(N)}\|_1^2+\|{\bf B}^{(N)}\|_1^2\big)^{5/2}\Big).\label{d2b}
\end{align}

We observe that the r.h.s.~of \rf{d2b} involves
powers of the sum $\|{\bf V}^{(N)}\|_1^2+\|{\bf B}^{(N)}\|_1^2$ that are too
high (larger than 1) to guarantee the time integrability
of the r.h.s.~(apparently, this also happens for any larger $s$).
Thus, the quadratic nonlinearity in \rf{vb}
prevents us from demonstrating, by using \rf{d2b} directly, that the derivatives
$\d{\bf V}_{\bf n}^{(N)}\!/\d t$ and $\d{\bf B}_{\bf n}^{(N)}\!/\d t$ are
uniformly (over the resolution parameter $N$) equicontinuous on $[0,T]$
for a given $T>0$.
Nevertheless, a subtler reasoning gives an opportunity to establish
the desired result, using the bound \rf{d2b} for times, when
$\|{\bf V}\|_1+\|{\bf B}\|_1$ is finite. We show this
in section~\ref{eqc}.

\subsection{The singularity set of solutions to equations
of magnetohydrodynamics}\label{str}

It was established in \cite{Le,FGT} that there exists an open set of times
such that the $H_m(\T)$ norms of weak solutions to the Navier--Stokes equation
are finite and continuous for all $m\ge1$, and
the complement has the Lebesgue measure zero, provided the initial condition
belongs to $H_1(\T)$. We apply now the approach of \cite{FGT}
to the equations of magnetohydrodynamics \rf{vb}.

We have shown that if $\|{\bf V}\|_s+\|{\bf B}\|_s$
is finite for some $s>1/2$ at a certain time $t=t_0$, then for
$t_0<t<t_0+t_*(t_0)$ (see \rf{tx}) the solution consists of space-analytic
vector fields.

{\it Definition}.
For $s>1/2$, an open time interval $t_0<t<t_1$ such that $0\le t_0<t_1\le T$
is called an {\it $H_s$-regularity interval} for a solution to \rf{vb}, if on
this interval ${\bf V(x},t)$ and ${\bf B(x},t)$ belong to $\dot H_s(\T)$ and
depend continuously on time in the norm $\|\cdot\|_s$. The open interval is
called a {\it maximal $H_s$-regularity interval}, if no larger $H_s$-regularity
interval including $(t_0,t_1)$ exists in $[0,T]$ for this solution.

{\it Definition}.
Suppose $\bf V$ and $\bf B$ belong to $\dot H_s(\T)$ at a time
$t=t_0$ for some $s>1/2$. The open time interval $t_0<t<t_0+t_*(t_0)$
is called the {\it time interval of guaranteed space analyticity}. An open
interval is called a {\it maximal interval of space analyticity}, if there
does not exist in $[0,T]$ any larger open interval on which this solution
is space-analytic at any time, including $(t_0,t_1)$.

{\it Theorem 6.}
Suppose ${\bf V}^{\rm(init)}$ and ${\bf B}^{\rm(init)}$ belong to $\dot H_s(\T)$
for some $s>1/2$. We focus on the solution for $t\le T$. Let $\O$ be
the intersection of $(0,T)$ with the union of all intervals of guaranteed space
analyticity of the solution, such that their left ends $t$ satisfy $0\le t<T$.

$i.$ The set $\O$ is open. The Lebesgue measure of the complement
$[0,T]\backslash\O$ is zero.

$ii.$ For any $p>1/2$, maximal $H_p$-regularity intervals coincide with maximal
intervals of space analyticity.

$iii.$ Each maximal $H_p$-regularity interval is also a maximal interval
of $H_{3/2}$-regularity of the transformed solutions $\wv$ and $\wb$
(see \rf{svb} and \rf{wvb}) to the auxiliary problem \rf{tr}.

{\it Proof}. By virtue of the energy inequality \rf{en}, the set
$$\S=\{t\in[0,T]\,\big|\,\|{\bf V}\|_s+\|{\bf B}\|_s
=\infty~\mbox{for all}~s>1/2\}$$
has the Lebesgue measure zero. Any time $t$ in the complement
$[0,T]\backslash\S$ can serve as the left end of an open interval $\O(t)$
of guaranteed space analyticity (by \rf{qu}, $\O(t)$ has an empty intersection
with~$\S$). The union of open intervals
$\O=\cup_{t\in[0,T]\backslash\S}(\O(t)\cap(0,T))$ is open.
Any connected component of $\O$ is a maximal interval of space
analyticity of the solution. The set $([0,T]\backslash\S)\backslash\O$ consists
of end points of such intervals, and hence it is at most countable (since each
interval
contains a rational point) and has the Lebesgue measure zero. This proves $i$.

Let us consider a maximal $H_p$-regularity interval $(l,r)$ for $p>1/2$.
Any point in $(l,r)$ is the left end of an interval of guaranteed space
analyticity. These intervals cover the entire interval $(l,r)$. If, otherwise,
$\hat t\in(l,r)$ is not covered, then $t_*(t_k)<\hat t-t_k$ for any
monotonically increasing sequence $t_k\to\hat t$, and a contradiction arises:
by the time continuity of $\|{\bf V}\|_p$ and $\|{\bf B}\|_p$ on $(l,r)$,
the norms have a uniform upper bound in a sufficiently
short closed interval $[\hat t-\epsilon,\hat t]$ and thus the lengths
$t_*(t)$ of the intervals of guaranteed space analyticity with the left ends
$t\in[\hat t-\epsilon,\hat t]$ have a uniform positive bound from below. Thus,
$(l,r)$ belongs to a maximal interval of space analyticity.

To prove the converse, we note that at each point of a maximal interval
of space analyticity, which we now denote $(l,r)$, $\|{\bf V}\|_p$ and
$\|{\bf B}\|_p$ are finite for any $p>1/2$ and hence the solution belongs to
$H_p(\T)$. Thus, to establish that $(l,r)$ belongs to a maximal $H_p$-regularity
interval, it suffices to show that the norms are continuous in time.

To do this, we first show that $\|{\bf V}\|_p$ and $\|{\bf B}\|_p$ are uniformly
bounded for any $p>0$ on a closed subinterval $[l+\epsilon,r-\epsilon]$, where
$\epsilon>0$ is sufficiently small. By construction, the maximal interval $(l,r)$
is covered by open intervals $\O(t)$ of guaranteed space analyticity. Hence,
we can choose a finite coverage $\{\O(t_k)\,|\,1\le k\le K\}$ of
$[l+\epsilon,r-\epsilon]$. The function
$$\widetilde n_p(t-t_k;t_k)=\left\{\begin{array}{ll}
1/n_p(t-t_k;t_k),~~~~&t_k<t<t_k+t_*(t_k),\\
0,&t\le t_k,\mbox{~~or~~}t\ge t_k+t_*(t_k)\end{array}\right.$$
is continuous on $\mathbb R$ (see \rf{qu}). Consequently,
$\max_{\,1\le k\le K}\widetilde n_p(t-t_k;t_k)$ is also continuous
on $\mathbb R$ and hence it admits its minimum on the closed interval
$[l+\epsilon,r-\epsilon]$. The minimum is strictly positive, since
its vanishing at a certain $t$ would indicate that this $t$ is outside of each
of the $K$ intervals $\O(t_k)$ covering the subinterval. Therefore,
$\|{\bf V}\|_p^2+\|{\bf B}\|_p^2\le1/\max_{\,1\le k\le K}\widetilde n_p(t-t_k;t_k)$
is uniformly bounded on $[l+\epsilon,r-\epsilon]$.

Second, for any $s$ we establish the time continuity of the solution
in the norm $\|\cdot\|_s$ on the same closed subinterval.
Due to convergence of the Fourier harmonics $\hV_{\bf n}^{(N_j)}$ and
$\hB_{\bf n}^{(N_j)}$ on $[0,T]$ when $N_j\to\infty$, \rf{UB} implies
$$|\hV_{\bf n}(t_1)-\hV_{\bf n}(t_2)|\le\widehat C|{\bf n}|^2|t_1-t_2|,\qquad
|\hB_{\bf n}(t_1)-\hB_{\bf n}(t_2)|\le\widehat C|{\bf n}|^2|t_1-t_2|.$$
Thus,
\ba&\|\hV(t_1)-\hV(t_2)\|_s^2+\|\hB(t_1)-\hB(t_2)\|_s^2
=\!\sum_{{\bf n}\ne0}\!\Big(\,\big|\hV_{\bf n}(t_1)-\hV_{\bf n}(t_2)\big|^2
+\big|\hB_{\bf n}(t_1)-\hB_{\bf n}(t_2)\big|^2\Big)|{\bf n}|^{2s}\\
&\le|t_1-t_2|\widehat C\sum_{{\bf n}\ne0}\!\Big(\,
\big|\hV_{\bf n}(t_1)\big|+\big|\hV_{\bf n}(t_2)\big|
+\big|\hB_{\bf n}(t_1)\big|+\big|\hB_{\bf n}(t_2)\big|\Big)|{\bf n}|^{2s+2}\\
&\le|t_1-t_2|\widehat C\Big(8\sum_{{\bf n}\ne0}\!|{\bf n}|^{-4}\Big)^{\!\!1/2}
\max_{l+\epsilon\le t\le r-\epsilon}\Big(\sum_{{\bf n}\ne0}\!\big(\big|\hV_{\bf n}(t)\big|^2
+\big|\hB_{\bf n}(t)\big|^2\big)|{\bf n}|^{4s+8}\Big)^{\!\!1/2}\\
&\le|t_1-t_2|\widehat C\Big(8\sum_{{\bf n}\ne0}\!|{\bf n}|^{-4}\Big)^{\!\!1/2}
\max_{l+\epsilon\le t\le r-\epsilon}(\|\hV(t)\|_{2s+4}+\|\hB(t)\|_{2s+4}),\end{align*}
which proves the continuity on $[l+\epsilon,r-\epsilon]$, since
$\|\hV\|_{2s+4}$ and $\|\hB\|_{2s+4}$ are uniformly bounded on this closed
interval. Since $\epsilon>0$ is arbitrary, the solutions are continuous
in time in the norm $\|\cdot\|_s$ on the entire maximal interval of space
analyticity.

If for some $s>1/2$ the sum $\|\hV(t_k)\|_s\!+\!\|\hB(t_k)\|_s$ is bounded
for a sequence of $t_k\to r$, then by \rf{tx} the intervals of guaranteed space
analyticity beginning at $t_k$ have lengths bounded from below by a positive
constant. This contradicts with the assumption that $(l,r)$ is a maximal
interval of space analyticity of the solutions. Therefore, in every such interval
$\lim_{\,t\to r}\|\hV(t)\|_s+\|\hB(t)\|_s=\infty$. Statement $ii$ is proven.

The transformation \rf{svb} of the Fourier coefficients introduced
in section~\ref{bo} can be implemented for any $\delta>0$, provided the Fourier
series is an analytic function that has a strictly positive size of the region
of analyticity. In particular, such a transformation and construction
of the transformed series $\wv$ and $\wb$ \rf{wvb} is possible everywhere
in $\O$, the resultant fields $\wv$ and $\wb$ belonging to $\dot H_{3/2}(\T)$.
Consequently, any connected component of $\O$ is also a maximal interval
of $H_{3/2}$-regularity of the solutions $\wv$ and $\wb$ to the auxiliary
problem \rf{tr}. The proof of the Theorem is completed. $\RHD$

\subsection{Application of \rf{d2b} for proving the bounds
\rf{dBV} for weak solutions}\label{eqc}

We focus on the subsequence of the Fourier--Galerkin approximants
${\bf V}^{(N_j)},\,{\bf B}^{(N_j)}$, whose limit is the weak solution
at hand to \rf{vb} on a certain time interval $0\le t\le T$.
We prove here equicontinuity of the time derivatives
$\d{\bf V}_{\bf n}^{(N_j)}\!/\d t$ and
$\d{\bf B}_{\bf n}^{(N_j)}\!/\d t$ for each wave vector~$\bf n$ on any closed
subinterval of a maximal interval of space analyticity. To carry over
the bounds \rf{dBV} to weak solutions, we apply a technical Theorem 7.

{\it Theorem 7.} Let $l<t<r$ be a maximal interval of space analyticity and
$\epsilon$ an arbitrary number satisfying $0<\epsilon<(r-l)/3$.
The Fourier--Galerkin
approximants ${\bf V}^{(N_j)}$, ${\bf B}^{(N_j)}$, that tend to the weak
solution to \rf{vb} under consideration, converge in $H_1(\T)$ uniformly
on the time interval $l+2\epsilon\le t\le r-\epsilon$, and thus
$\|{\bf V}^{(N_j)}\|_1$ and $\|{\bf B}^{(N_j)}\|_1$ are uniformly
(over $N_j$) bounded on this interval.

{\it Proof.}
Let us consider a closed subinterval $l+\epsilon\le t\le r-\epsilon$
of a maximal interval $l<t<r$ of space analyticity, where $0\le l<r\le T$.
We exploit compactness of the embedding $H_2(\T)\subset H_1(\T)$.
By Theorem~6,
\BE m_s=\max_{l+\epsilon\le t\le r-\epsilon}(\|{\bf V}\|_s^2+\|{\bf B}\|_s^2)
<\infty.\EE{mx}
Since the maximum $m_2$ is finite,
\se{uob}\BE\sum_{|{\bf n}|>(m_2/\zeta)^{1/2}}|{\bf n}|^2(|\hV_{\bf n}|^2
+|\hB_{\bf n}|^2)<\zeta\quad\mbox{for all $t$
in the subinterval~~}l+\epsilon\le t\le r-\epsilon,\EE{bT}
where $\zeta>0$ is arbitrary. The numbers
$$M_1^{(N)}=\int_{l+\epsilon}^{l+2\epsilon}\!\!\big(
\|{\bf V}^{(N)}\|_2^2+\|{\bf B}^{(N)}\|_2^2\big)^{1/3}\d t$$
do not exceed the r.h.s.~of \rf{qq} for $s=2$, which is independent of $N$.
Thus, for each $N$, there exists a point $t=\tau^{(N)}$ in the interval
$l+\epsilon\le t\le l+2\epsilon$, at which
$$\|{\bf V}^{(N)}\|_2^2+\|{\bf B}^{(N)}\|_2^2\le
M_2=\sup_N\,(M^{(N)}_1\!/\epsilon)^3,$$
whereby
\BE\sum_{|{\bf n}|>(M_2/\zeta)^{1/2}}|{\bf n}|^2(|\hV^{(N)}_{\bf n}|^2
+|\hB^{(N)}_{\bf n}|^2)\le\zeta\quad\mbox{at~~}t=\tau^{(N)}.\EE{bNT}
Due to the weak convergence ${\bf V}^{(N_j)}\to\bf V$ and ${\bf B}^{(N_j)}\to\bf B$
for $N_j\to\infty$, there exists $M_3$ such that
\BE\sum_{|{\bf n}|\le(\max(m_2,M_2)/\zeta)^{1/2}}\!\!|{\bf n}|^2
(|\hV_{\bf n}-\hV^{(N_j)}_{\bf n}|^2+|\hB_{\bf n}-\hB^{(N_j)}_{\bf n}|^2)\le\zeta
\quad\mbox{everywhere on~}[0,T]\EE{bD}\end{subequations}
for all $N_j\ge M_3$. Together, inequalities \rf{uob} imply that, given $\zeta>0$
and $\epsilon>0$, we can find for any $N_j\ge M_3$ a point $t=\tau^{(N_j)}$
in the interval $l+\epsilon\le t\le l+2\epsilon$, at which
the norms of the discrepancies ${\bf u}={\bf V}-{\bf V}^{(N_j)}$,
${\bf a}={\bf B}-{\bf B}^{(N_j)}$
are controlled:
\BE\|{\bf u}\|_1^2+\|{\bf a}\|_1^2\le3\zeta.\EE{at}

It is convenient to split the discrepancies in two parts:
$${\bf u}=\sum_{\bf n}\hu_{\bf n}\e^{\i\bf n\cdot x}
={\bf u}^<+{\bf u}^>,\qquad{\bf u}^<=\!\!\sum_{|{\bf n}|\le N_j}
(\hV_{\bf n}-\hV_{\bf n}^{(N_j)})\e^{\i\bf n\cdot x},\qquad
{\bf u}^>=\!\!\sum_{|{\bf n}|>N_j}\hV_{\bf n}\e^{\i\bf n\cdot x};$$
$${\bf a}=\sum_{\bf n}\ha_{\bf n}\e^{\i\bf n\cdot x}
={\bf a}^<+{\bf a}^>,\qquad{\bf a}^<=\!\!\sum_{|{\bf n}|\le N_j}
(\hB_{\bf n}-\hB_{\bf n}^{(N_j)})\e^{\i\bf n\cdot x},\qquad
{\bf a}^>=\!\!\sum_{|{\bf n}|>N_j}\hB_{\bf n}\e^{\i\bf n\cdot x}.$$
Fourier coefficients of ${\bf u}^<$ and ${\bf a}^<$
satisfy the following equations for $|{\bf n}|\le N_j$:
\se{dis}\be{\d\hu_{\bf n}\over\d t}=&-\nu|{\bf n}|^2\hu_{\bf n}
-\i\sum_{\bf k}\P\big((\hu_{\bf n-k}\!\cdot\!{\bf k})\hV_{\bf k}
+(\hV_{\bf n-k}\!\cdot\!{\bf k})\hu_{\bf k}
-(\hu_{\bf n-k}\!\cdot\!{\bf k})\hu_{\bf k}\nonumber\\
&\hspace*{9em}-(\ha_{\bf n-k}\!\cdot\!{\bf k})\hB_{\bf k}
-(\hB_{\bf n-k}\!\cdot\!{\bf k})\ha_{\bf k}
+(\ha_{\bf n-k}\!\cdot\!{\bf k})\ha_{\bf k}\big),\label{hV}\\
{\d\ha_{\bf n}\over\d t}=&-\eta|{\bf n}|^2\ha_{\bf n}
+\i{\bf n}\times\sum_{\bf k}\,\big(\hu_{\bf k}\!\times\hB_{\bf n-k}
+\hV_{\bf k}\!\times\ha_{\bf n-k}-\hu_{\bf k}\!\times\ha_{\bf n-k}\big).
\label{hB}\end{align}\end{subequations}

By \rf{bT}, for the fixed $\zeta$ and $\epsilon$, discrepancies for
$N_j>\max(M_3,(m_2/\zeta)^{1/2})$ satisfy
\BE\|{\bf u}^>\|_1^2+\|{\bf a}^>\|_1^2<\zeta\quad\mbox{for all $t$
in the subinterval~~}l+\epsilon\le t\le r-\epsilon.\EE{xt}
Scalar multiplying \rf{hV} and \rf{hB} for $|{\bf n}|\le N_j$
by $|{\bf n}|^2\,\hu_{\bf-n}$ and $|{\bf n}|^2\,\ha_{\bf-n}$, respectively,
and summing up the results yields
\se{li}\be&{1\over2}\,{\d\over\d t}(\|{\bf u}^<\|_1^2+\|{\bf a}^<\|_1^2)
+\nu\|{\bf u}^<\|_2^2+\eta\|{\bf a}^<\|_2^2\nonumber\\
&\hspace*{1em}\le C'_{1/2}\Big(\!\|{\bf u}^<\|_2\big(
\|{\bf u}\|_1\|{\bf V}\|_{3/2}\!+\|{\bf a}\|_1\|{\bf B}\|_{3/2}\big)
+\|{\bf a}^<\|_2\big(\|{\bf u}\|_1\|{\bf B}\|_{3/2}
\!+\|{\bf a}\|_1\|{\bf V}\|_{3/2}\big)\!\Big)\label{z1}\\
&\hspace*{1em}+C'_{1/2}\Big(\!\|{\bf u}^<\|_2\big(
\|{\bf u}\|_{3/2}(\|{\bf V}\|_1\!+\|{\bf u}\|_1)
+\|{\bf a}\|_{3/2}(\|{\bf B}\|_1\!+\|{\bf a}\|_1)\big)\nonumber\\
&\hspace*{3.5em}+\|{\bf a}^<\|_2\big(
\|{\bf u}\|_{3/2}(\|{\bf B}\|_1\!+\|{\bf a}\|_1)
+\|{\bf a}\|_{3/2}(\|{\bf V}\|_1\!+\|{\bf u}\|_1)\big)\!\Big).
\label{z2}\end{align}\end{subequations}
We bound the two sums in the r.h.s.~of this inequality
on the subinterval $l+\epsilon\le t\le r-\epsilon$ separately.
The first one, \rf{z1}, by Young's inequality, does not exceed
\se{xx}\be&(\nu\|{\bf u}^<\|_2^2+\eta\|{\bf a}^<\|_2^2)/4
+(C'_{1/2})^2(1/\nu+1/\eta)(\|{\bf u}\|_1^2
+\|{\bf a}\|_1^2)(\|{\bf V}\|_{3/2}^2\!+\|{\bf B}\|_{3/2}^2)\nonumber\\
\le\,&(\nu\|{\bf u}^<\|_2^2+\eta\|{\bf a}^<\|_2^2)/4
+(C'_{1/2})^2\,m_{3/2}(1/\nu+1/\eta)
(\|{\bf u}^<\|_1^2+\|{\bf a}^<\|_1^2+\zeta)\label{x1}\end{align}
(\rf{mx} for $s=3/2$ and \rf{xt} have been used).
The second sum, \rf{z2}, has an upper bound
\be&C'_{1/2}(\|{\bf u}^<\|_2+\|{\bf a}^<\|_2)
\Big(\!\big(\|{\bf u}\|_{3/2}^2+\|{\bf a}\|_{3/2}^2\big)
\big((\|{\bf V}\|_1\!+\|{\bf u}\|_1)^2+(\|{\bf B}\|_1\!+\|{\bf a}\|_1)^2
\big)\!\Big)^{\!\!1/2}\nonumber\\
\le\,&C'_{1/2}(\|{\bf u}^<\|_2+\|{\bf a}^<\|_2)\Big(\!2\big(\|{\bf u}^<\|_2
\|{\bf u}^<\|_1\!+\|{\bf a}^<\|_2\|{\bf a}^<\|_1\!+\sqrt{m_2\zeta}\,\big)
\big(\|{\bf u}^<\|_1^2\!+\|{\bf a}^<\|_1^2\!+\zeta\!+m_1\big)\!\Big)^{\!\!1/2}\nonumber\\
\le\,&2C'_{1/2}\big(\|{\bf u}^<\|_2^2+\|{\bf a}^<\|_2^2\big)^{\!1/2}
\Big(\!\big(\|{\bf u}^<\|_2^2+\|{\bf a}^<\|_2^2\big)^{\!1/4}
\big(\|{\bf u}^<\|_1^2+\|{\bf a}^<\|_1^2\big)^{\!1/4}+(m_2\zeta)^{1/4}\Big)
\nonumber\\
&\hspace*{2em}\times\big(\|{\bf u}^<\|_1^2\!
+\|{\bf a}^<\|_1^2\!+\zeta\!+m_1\big)^{1/2}\nonumber\\
\le\,&{\min(\nu,\eta)\over2}\big(\|{\bf u}^<\|_2^2\!+\|{\bf a}^<\|_2^2\big)
+{108(C'_{1/2})^4\over(\min(\nu,\eta))^3}
\big(\|{\bf u}^<\|_1^2+\|{\bf a}^<\|_1^2\big)
\big(\|{\bf u}^<\|_1^2+\|{\bf a}^<\|_1^2\!+\zeta\!+m_1\big)^2\nonumber\\
&\hspace*{2em}+{4(C'_{1/2})^2(m_2\zeta)^{1/2}\over\min(\nu,\eta)}\big(
\|{\bf u}^<\|_1^2+\|{\bf a}^<\|_1^2+\zeta\!+m_1\big)
\label{x2}\end{align}\end{subequations}
(Young's inequality has been again employed).

The two bounds \rf{xx} give rise to a differential inequality for
$\psi=\|{\bf u}^<\|_1^2+\|{\bf a}^<\|_1^2$ of the form
\BE{\d\psi\over\d t}\le A_3\psi^3+A_2\psi^2+A_1\psi+A_0.\EE{DI}
The constant $A_0$ in \rf{DI} is proportional to $\zeta^{1/2}$, and, for small
$\zeta$, the three remaining constants $A_i$ are O(1). Initial conditions
\rf{at} at $t=\tau^{(N_j)}$ are O($\zeta$). Thus, \rf{DI} implies an
O($\zeta^{1/2}$) upper bound for $\psi$ on the O(1)-long time
interval $\tau^{(N_j)}\le t\le r-\epsilon$. Consequently, $\psi\to0$ uniformly
on the interval $l+2\epsilon\le t\le r-\epsilon$ in the limit $N_j\to\infty$
and $\zeta\to0$, and hence
$\|{\bf u}\|_1^2+\|{\bf a}\|_1^2=\psi+\|{\bf u}^>\|_1^2+\|{\bf a}^>\|_1^2\to0$.
This proves the Theorem.~~$\RHD$

We are now in a position to achieve the goal of this section.

{\it Theorem 8.} Time derivatives of weak solutions to
the problem \rf{vb} obey the bounds \rf{dBV}.

{\it Proof.} Evidently,
\ba&\left|{\d\over\d t}\!\left(
\nu\left\|{\partial{\bf V}^{(N)}\over\partial t}\right\|_{s+1}^2
\!\!+\eta\left\|{\partial{\bf B}^{(N)}\over\partial t}\right\|_{s+1}^2\right)
\right|\le2\nu\left\|{\partial{\bf V}^{(N)}\over\partial t}\right\|_{s+2}
\left\|{\partial^2{\bf V}^{(N)}\over\partial t^2}\right\|_s
\!\!+2\eta\left\|{\partial{\bf B}^{(N)}\over\partial t}\right\|_{s+2}
\left\|{\partial^2{\bf B}^{(N)}\over\partial t^2}\right\|_s\nonumber\\
&\quad\le{1\over2}\left\|{\partial^2{\bf V}^{(N)}\over\partial t^2}\right\|_s^2
\!\!+{1\over2}\left\|{\partial^2{\bf B}^{(N)}\over\partial t^2}\right\|_s^2
\!\!+2\nu^2\left\|{\partial{\bf V}^{(N)}\over\partial t}\right\|_{s+2}^2
\!\!+2\eta^2\left\|{\partial{\bf B}^{(N)}\over\partial t}\right\|_{s+2}^2\!\!\!.
\end{align*}
Hence, \rf{d2b} for $s<-7/2$ and \rf{clo} imply
\ba\left\|{\partial^2{\bf V}^{(N_j)}\over\partial t^2}\right\|_s^2
\!\!+\left\|{\partial^2{\bf B}^{(N_j)}\over\partial t^2}\right\|_s^2
\le&\,4\max(\nu^2,\eta^2)\,
\widetilde C''\big(\|{\bf V}^{(N_j)}\|_1^2+\|{\bf B}^{(N_j)}\|_1^2
+(\|{\bf V}^{(N_j)}\|_1^2+\|{\bf B}^{(N_j)}\|_1^2)^{3/2}\big)\\
&+40c^2_{-s-2}\,
\widetilde C''\big((\|{\bf V}^{(N_j)}\|_1^2+\|{\bf B}^{(N_j)}\|_1^2)^2
+(\|{\bf V}^{(N_j)}\|_1^2+\|{\bf B}^{(N_j)}\|_1^2)^{5/2}\big)\\
\le&M',\end{align*}
where $M'$ denotes the finite on the interval $l+2\epsilon\le t\le r-\epsilon$
supremum (over $N_j$) of the middle part of this inequality. Therefore,
$$\Big|{\d\hV_{\bf n}^{(N_j)}\over\d t}(t')-{\d\hV_{\bf n}^{(N_j)}\over\d t}(t'')
\Big|\le|{\bf n}|^{-s}\sqrt{M'}|t'-t''|,\qquad
\Big|{\d\hB_{\bf n}^{(N_j)}\over\d t}(t')-{\d\hB_{\bf n}^{(N_j)}\over\d t}(t'')
\Big|\le|{\bf n}|^{-s}\sqrt{M'}|t'-t''|$$
for any $t'$ and $t''$ belonging to this interval, whereby, for each wave
vector~$\bf n$, $\d\hV_{\bf n}^{(N_j)}\!/\d t$ and
$\d\hB_{\bf n}^{(N_j)\!}/\d t$ are
equicontinuous on this time interval. Relying on the Arzel\`a--Ascoli theorem
and employing the diagonal process, we construct a subsequence
$N_j\to\infty$ such that, for any $\bf n$, the derivatives
$\d\hV_{\bf n}^{(N_j)}\!/\d t$ and $\d\hB_{\bf n}^{(N_j)}\!/\d t$ converge
uniformly on the interval $l+2\epsilon\le t\le r-\epsilon$ to some
continuous limit functions $\bm\phi_{\bf n}^{\bf V}(t)$
and $\bm\phi_{\bf n}^{\bf B}(t)$, respectively. Taking the limit $N_j\to\infty$
in the identities
$$\int_{l+2\epsilon}^\tau{\d\hV_{\bf n}^{(N_j)}\over\d t}\d t=
\hV_{\bf n}^{(N_j)}(\tau)-\hV_{\bf n}^{(N_j)}(l+2\epsilon),\qquad
\int_{l+2\epsilon}^\tau{\d\hB_{\bf n}^{(N_j)}\over\d t}\d t=
\hB_{\bf n}^{(N_j)}(\tau)-\hB_{\bf n}^{(N_j)}(l+2\epsilon),$$
we obtain relations
$$\int_{l+2\epsilon}^\tau\bm\phi_{\bf n}^{\bf V}\d t=
\hV_{\bf n}(\tau)-\hV_{\bf n}(l+2\epsilon),\qquad
\int_{l+2\epsilon}^\tau\bm\phi_{\bf n}^{\bf B}\d t=
\hB_{\bf n}(\tau)-\hB_{\bf n}(l+2\epsilon),$$
equivalent to $\bm\phi_{\bf n}^{\bf V}=\d\hV_{\bf n}/\d t$ and
$\bm\phi_{\bf n}^{\bf B}=\d\hB_{\bf n}/\d t$.

Thus, for a subsequence of $N_j\to\infty$,
$\d\hV_{\bf n}^{(N_j)}\!/\d t$ and $\d\hB_{\bf n}^{(N_j)}\!/\d t$
converge to the derivatives of the harmonics $\d\hV_{\bf n}/\d t$ and
$\d\hB_{\bf n}/\d t$, respectively, on the interval
$l+2\epsilon\le t\le r-\epsilon$. Recalling that $\epsilon>0$ is
an arbitrary sufficiently small number and considering now the problem
for a sequence $\epsilon_k\to0$, we establish the convergence
$\d\hV_{\bf n}^{(N_j)}\!/\d t\to\d\hV_{\bf n}/\d t$ and
$\d\hB_{\bf n}^{(N_j)}\!/\d t\to\d\hB_{\bf n}/\d t$ on the entire
$H_1$-regularity interval $l<t<r$ for a subsequence of $N_j\to\infty$
(for which we keep the notation~$N_j$), employing the diagonal process
on increasing the interval. Since the $H_1$-regularity intervals are countable,
employing again the diagonal process, we can distill a subsequence, for which
the convergence occurs on the entire set $\O$, i.e., almost everywhere
in the interval $0\le t\le T$.

To show that the a priori bounds \rf{dBV} hold true for the weak solutions
of the problem \rf{vb}, we note that the Fourier series for approximants
truncated at a level $M$ satisfy \rf{dBV}, e.g., \rf{ds} implies
$$\sum_{|{\bf n}|\le M}|{\bf n}|^{2s}\big(\big|\d\hV_{\bf n}^{(N_j)}\!/\d t\big|^2
+\big|\d\hB_{\bf n}^{(N_j)}\!/\d t\big|^2\big)\le D_s^{(3)}.$$
Convergence of the time derivatives of individual harmonics almost
everywhere being proven, this inequality, for a fixed $M$, holds upon
taking the limit $N_j\to\infty$, and then the inequality \rf{ds}
for the weak solutions follows almost everywhere since $M$ is arbitrary.
Similarly, \rf{das} implies
$$\int_{t_0}^T\!\!\Big(\sum_{|{\bf n}|\le M}|{\bf n}|^{2s}\big(\big|\d\hV^{(N_j)}_{\bf n}\!/\d t\big|^2
+\big|\d\hB^{(N_j)}_{\bf n}\!/\d t\big|^2\big)\!\Big)^{\!\alpha_{s+2}/2}\,\d t\le D_s^{(1)}.$$
For a fixed $M$, the sum in the integrand converges almost
everywhere for $t\!\le\!T$ to the analogous sum for the Fourier coefficients
of the weak solution. Hence, by Fatou's lemma (see, e.g., \cite{We}),
the inequality holds true in the limit $N_j\to\infty$. Truncations
being arbitrary, \rf{das} follows for weak solutions, as desired, when
$M\to\infty$. The inequality \rf{dgs} is proven by a similar argument.

Finally, we use a similar approach to demonstrate the Wiener norm bound
\rf{dtm} for weak solutions. We have proven a stronger a priori bound
$$\int_{t_0}^T\!\!\Big(\sum_{\bf n}|{\bf n}|^s\big(\big|\d\hV^{(N_j)}_{\bf n}\!/\d t\big|
+\big|\d\hB^{(N_j)}_{\bf n}\!/\d t\big|\big)\!\Big)^{\!\alpha_{s+7/2}}\,\d t\le D_s^{(4)}$$
implying \rf{dtm}. In view of the convergence
$\d\hV^{(N_j)}_{\bf n}\!/\d t\to\d\hV_{\bf n}/\d t$
and $\d\hB^{(N_j)}_{\bf n}\!/\d t\to\d\hB_{\bf n}/\d t$ at almost all times
when $N_j\to\infty$, by Fatou's lemma this inequality holds true for truncated
sums for the time derivatives of the Fourier coefficients
of weak solutions. Letting the truncation parameter tend to infinity
proves \rf{dtm} for the weak solution.~~$\RHD$

\section{Concluding remarks}\label{cn}

The similarity of the quadratic nonlinearity of the terms describing advection
and the Lorentz force in the Navier--Stokes equation and in the magnetic
induction equation has enabled us to carry over the
results of the theory of the Navier--Stokes equation to the system of equations
of magnetohydrodynamics. Namely, applying the techniques of \cite{FT} we have
shown that the MHD solutions instantaneously acquire space analyticity,
provided initially they have a minimum regularity of $H_s(\T)$ for $s>1/2$ (see
section \ref{io}). Next, following \cite{Gev} we have introduced the auxiliary
problem \rf{tr} for vector fields, whose Fourier series involve transformed
coefficients (section \ref{tra}). Solutions to the auxiliary problem admit
the energy-like a priory bound \rf{bou} (section \ref{eb}) that yields
an integral bound for the $H_{3/2}(\T)$ norm of these
solutions. The inverse of this norm serves as a lower bound for the size
of the spatial analyticity region of the solutions $\bf V,B$ to the original
MHD problem; we thus obtain a simple proof that $\bf V,B$ are space-analytic
vector fields at almost all times. Relying on space analyticity, we derive
a priori bounds for $H_s(\T)$ norms of the solutions for arbitrary indices
$s$ (section \ref{bSo}), that are direct generalisations of the bounds derived
in \cite{FGT} in the hydrodynamic setup. An integral a priori bound
for the maximum of the flow velocity in the cube of periodicity was also
presented {\it ibid}. We have expanded this result by constructing bounds for
the Wiener algebra norms (i.e., the sums of absolute values of the Fourier
coefficients) of the fields $(-\nabla^2)^{s/2}\bf V$ and
$(-\nabla^2)^{s/2}\bf B$ for arbitrary $s>-1/2$ (section \ref{bWi}).
It is notable that three
independent approaches (the original one of \cite{FGT}, the one relying
on ladder inequalities \cite{G1,G2,G3,G4}, and the present one) yield the same
exponents $\alpha_m$ in \rf{NSb}, suggesting that these values are optimal and
cannot be improved unless construction of the bounds is based on new,
significantly different ideas. Finally, we have derived similar a priori
integral bounds for the Sobolev space and Wiener algebra norms (section
\ref{dt}) of $\d{\bf V}\!/\d t$ and $\d{\bf B}/\d t$. Proving that the a priori
bounds hold for the time derivatives of the weak MHD solutions
(section~\ref{bw}) is considerably more involved than those for the solution
itself. This has required identifying the structure of the singularity set
of the solution in the time domain and proving convergence in the $H_1(\T)$
norm of the relevant subsequence of the Fourier--Galerkin approximants
at times in the complement to this set (section \ref{str}). We have thus
demonstrated that the bounds for the Sobolev and Wiener norms of the MHD
solutions and their time derivatives stem from their space analyticity.

According to the present paradigm, the action of viscosity and diffusivity
hampers development of small-scale structures generated by the nonlinearity.
Thus, bounding the diffusive and nonlinear terms jointly may be expected to
result in more accurate bounds for a larger ``number of derivatives'' (i.e.,
for a higher-index Sobolev space norm). We have not
achieved this when estimating the time derivatives of the MHD solutions:
our bounds for the nonlinear advective terms are for the same index norms, as
for the dissipative terms. Indirectly this confirms that cancellation may be
possible with the sum residing in a higher-index Sobolev space,
our estimations then being too conservative.

The singularity set of a weak solution is the zero-measure complement
to the union $\O$ of its maximal intervals of space analyticity, or the union
of maximal $H_p$-regularity intervals for any $p>1/2$. If for a certain initial
condition weak MHD solutions are non-unique, their branching occurs only
at times belonging to the singularity set. It is unclear, whether any specific
techniques for constructing weak solutions favour some of them that are
in some sense ``better''. We may mention the following difference in
construction of weak solutions using their Galerkin approximations (as we have
done in this paper), or regularising the original system of MHD equations
\rf{vb}. Regularisation can be achieved by introducing the hyperdiffusivity
terms $-\varepsilon(-\nabla^2)^p\bf V$ and $-\varepsilon(-\nabla^2)^p\bf B$
into the r.h.s.~of \rf{ve} and \rf{be}, respectively, for $\varepsilon>0$ and
$p\ge5/4$ (see \cite{Li}). Like in the hydrodynamic setup, it is easy to show that
the regularised solutions ${\bf V_\varepsilon(x},t)$, ${\bf B_\varepsilon(x},t)$
are strong and unique, they depend continuously on $\varepsilon$, and any
sequence $\varepsilon_j\to0$ contains a subsequence, for which the regularised
solutions weakly converge to a weak solution to the original problem \rf{vb}.
Either such a limit weak solution is unique (i.e., a weak limit exists for
$\varepsilon\to0$), or a continuum of weak solutions exist for the initial
condition at hand. (This stems from the fact that $\varepsilon$ is not a
discrete parameter: If the limit is non-unique, then
there exist vector fields $\bf f_V$ and $\bf f_B$ such that
$$w(\varepsilon)={\rm Re}\int_{\T}({\bf V}_{\varepsilon}({\bf x},t)\cdot{\bf f_V}+
{\bf B}_{\varepsilon}({\bf x},t)\cdot{\bf f_B})\d{\bf x}$$
tends for some $t>0$
to distinct limits $w(\varepsilon_{k,j})\to w_k,~w_1<w_2$ for two sequences
$\varepsilon_{k,j}\to0$, $j\to\infty$, $k=1,2$ (see Figure 1). By suitably
rarefying the two sequences, we can render them intermittent:
$\varepsilon_{1,j}<\varepsilon_{2,j}<\varepsilon_{1,j+1}<\varepsilon_{2,j+1}$
for all $j$. Let $w$ satisfy $w_1+\gamma<w<w_2-\gamma$ for a sufficiently small
$\gamma>0$. Because of the weak convergence,
$|w(\varepsilon_{k,j})-w_k|<\gamma$ for a sufficiently large $J$ and all $j>J$
for both sequences $\varepsilon_{k,j}$. Continuity in $\varepsilon>0$ implies
that there exists a sequence $\varepsilon_{3,j}\to0$ such that
$\varepsilon_{1,j}<\varepsilon_{3,j}<\varepsilon_{2,j}$ and
\BE w(\varepsilon_{3,j})=w.\EE{ws}
There exists a subsequence of $\varepsilon_{3,j}$ for which the regularised
solutions converge to a weak solution, such that, evidently, \rf{ws} holds.
Since the open interval $(w_1,w_2)$ consists of a continuum of such $w$,
a continuum of weak solutions exist for the initial
condition at hand.) While simple changes in the proofs of Theorems 2 and 3
suffice to specialise them for the solutions obtained by the hyperdiffusive
regularisation of the problem \rf{vb}, our proof cannot be modified
straightforwardly to justify the bounds of Theorem 4 for time
derivatives of these weak solutions.

\pagebreak
\begin{figure}[t]
\centerline{\includegraphics[width=3in]{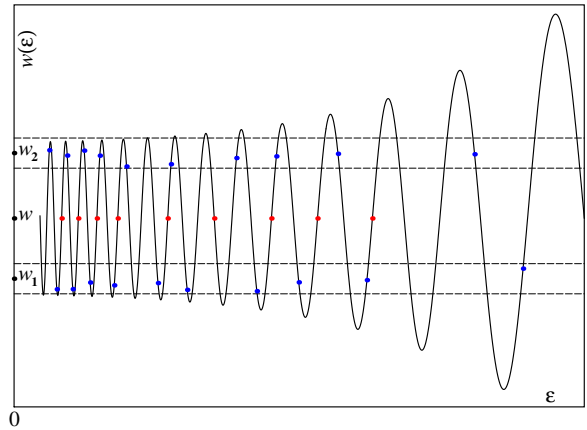}}
\caption{A sketch of behaviour of $w(\varepsilon)$ when two distinct limit weak
solutions to the MHD equations coexist (see
the text). Points $(\varepsilon_{k,j},w(\varepsilon_{k,j}))$ (blue dots)
tend to $(0,w_k)$ (black dots) for $j\to\infty$, $k=1,2$. The sequence
$(\varepsilon_{3,j},w(\varepsilon_{3,j}))\to(0,w)$ (red dots), where
$w(\varepsilon_{3,j})=w$, exists due to continuity of $w(\varepsilon)$ in
$\varepsilon$ for $\varepsilon>0$. Dashed lines: boundaries of the regions
$|w-w_k|\le\gamma$.}
\end{figure}

\end{document}